\documentclass[10pt, reqno]{amsart}
\usepackage[T1]{fontenc}
\usepackage[dvipsnames]{xcolor}
\usepackage{graphicx}
\usepackage{amssymb}
\usepackage[normalem]{ulem}

\usepackage{geometry}
\usepackage{hyperref}
\usepackage{amsmath,amssymb,mathtools}
\usepackage{stmaryrd}
\usepackage{amsfonts}
\usepackage{cases}
\usepackage{amsthm}
\usepackage{booktabs,multirow}
\usepackage{chemformula}
\usepackage{csquotes}

\newtheorem{remark}{Remark}

\def\ie{\emph{i.e.\/}}
\newcommand{\Forall}{\forall\;}
\newcommand{\RR}{\mathbb{R}}
\newcommand{\NN}{\mathbb{N}}

\newcommand{\Id}{{\rm Id}}

\newcommand{\Ec}{\mathcal{E}}
\newcommand{\Ic}{\mathcal{I}}
\newcommand{\Hc}{\mathcal{H}}

\newcommand{\Nc}{\mathcal{N}}
\renewcommand{\d}{{\rm d}}

\renewcommand{\lq}{\leqslant}

\newcommand{\br}{\mathbf{r}}
\newcommand{\bR}{\mathbf{R}}
\newcommand{\PP}{\mathbb{P}}
\newcommand{\R}{\mathbb{R}}
\newcommand{\wt}{\widetilde}

\begin{document}

\title{On basis set optimisation in quantum chemistry}
\author{Eric Canc\`es}\address{CERMICS, Ecole des Ponts and Inria Paris, 6 \& 8 avenue Blaise Pascal, 77455 Marne-la-Vall\'ee, France, \texttt{eric.cances@enpc.fr}}
\author{Genevi\`eve Dusson}\address{Laboratoire de Math\'ematiques de Besan\c{c}on, UMR CNRS 6623, Universit\'e Bourgogne Franche-Comt\'e, 16 route de Gray, 25030 Besan\c{c}on, France, \texttt{genevieve.dusson@math.cnrs.fr}}
\author{Gaspard Kemlin}\address{CERMICS, Ecole des Ponts and Inria Paris, 6 \& 8 avenue Blaise Pascal, 77455 Marne-la-Vall\'ee, France, \texttt{gaspard.kemlin@enpc.fr}}
\author{Laurent Vidal}\address{CERMICS, Ecole des Ponts and Inria Paris, 6 \& 8 avenue Blaise Pascal, 77455 Marne-la-Vall\'ee, France,
  \texttt{laurent.vidal@enpc.fr}}

\begin{abstract}
  In this article, we propose general criteria to construct optimal atomic
  centered basis sets in quantum chemistry. We focus in particular on two
  criteria, one based on the ground-state one-body density matrix of the
  system and the other based on the ground-state energy. The performance of
  these two criteria are then numerically tested and compared on a
  parametrized eigenvalue problem, which corresponds to a one-dimensional toy
  version of the ground-state dissociation of a diatomic molecule.
\end{abstract}

\makeatletter
\setlength{\@fptop}{0pt plus 1fil}
\setlength{\@fpbot}{0pt plus 1fil}
\makeatother
\maketitle

\section{Introduction}

In quantum chemistry, a central problem is the computation of the electronic
ground-state (GS) of a given molecular system. For many-electron systems, it is
not possible to solve the $N$-body Schr\"odinger equations and most
calculations are thus based on variational (e.g. Hartree--Fock) or
non-variational (e.g. coupled cluster) approximations of the latter, or on
Kohn--Sham density functional theory (DFT). For all these models, the
continuous equations (e.g. a nonlinear elliptic eigenvalue problem in the
Hartree--Fock or Kohn--Sham settings) are discretized into a finite-dimensional
approximation space. Approximation spaces constructed from atomic orbitals
(AO) basis sets~\cite{Helgaker2014-py,Olsen2021-so} have many advantages and
are therefore the most common choice in the quantum chemistry community. An AO
basis set consists of a collection of functions
$\chi=(\chi_\mu^z)_{z \in {\rm CE}, \; 1 \lq \mu \lq n_z}$ where ${\rm CE}$ is
a set of atomic numbers (e.g. ${\rm CE}=\{ 1,\ldots,92 \}$ for the
natural chemical elements of the periodic table), $n_z$ a positive
integer depending on the electronic shell-structure of the chemical element
with atomic number $z$, and $\chi_\mu^z \in H^1(\R^3)$ a fast decaying
function centered at the origin called an atomic orbital. Consider an atomic
configuration $\omega$ consisting of $M$ nuclei with nuclear charges
$z_1,\dots,z_M$ (in atomic units) and positions $\bR_1,\dots,\bR_M$ in the
three dimensional physical space. If the AO basis set $\chi$ is chosen by the
user, the (spatial component of the) one-electron finite-dimensional space in
which the chosen electronic structure model of a molecular system with atomic
configuration $\omega$ is discretized is
$$
{\mathcal X}_\omega\coloneqq\mbox{span}(\chi_1^{z_1}(\cdot-\bR_1),\dots,\chi_{n_{z_1}}^{z_1}(\cdot-\bR_1),\dots,\chi_1^{z_M}(\cdot-\bR_M),\dots,\chi_{n_{z_M}}^{z_M}(\cdot-\bR_M)).
$$
The accuracy of the approximation therefore crucially depends on the quality
of the AO basis set. The main advantage of AO basis sets is that only a small
number of AO per atoms (typically a dozen) are necessary to obtain a
relatively accurate result on most quantities of interest. This is in sharp
contrast with standard discretization methods used in the simulation of
partial differential equations such as finite-element methods. To make
connection with discretization methods used in mechanical and electrical
engineering, AO basis set discretization methods can be considered as spectral
methods~\cite{canuto2007spectral}, and share common features with the modal synthesis
method~\cite[Chapter 7]{charpentier1996methode}, \cite{charpentier1996component}. A drawback of AO basis sets is that conditioning quickly blows
up when increasing the size of the basis by including polarization and diffuse
basis functions, a problem known as overcompleteness~\cite{Lowdin1970-oz}. The
numerical errors due to this large condition number can deteriorate the
accuracy of the computed solutions and/or significantly increase computational
times. AO basis sets can therefore not be systematically improved in a
straightforward way.

\medskip

In the early days, AOs were Slater functions~\cite{Slater1930-bf}, with exponential decay and a cusp at the origin. It was then realized by Boys~\cite{Boys1950-nj} in 1950 that it was much more efficient from a computational viewpoint to use Gaussian-type orbitals (GTO), that are linear combinations of polynomials times Gaussian functions. Indeed the multi-center overlap, kinetic and Coulomb integrals
\begin{align*}
  &\int_{\R^3} \chi_i^{z_a}(\br-\bR_a)\, \chi_j^{z_b}(\br-\bR_b) \, \d\br, \qquad \int_{\R^3} \nabla \chi_i^{z_a}(\br-\bR_a)\cdot \nabla\chi_j^{z_b}(\br-\bR_b) \, \d\br, \\
  &\int_{\R^3} \frac{\chi_i^{z_a}(\br-\bR_a)\, \chi_j^{z_b}(\br-\bR_b)}{|\br-\bR_k|} \, \d\br, \qquad
  \int_{\R^3 \times \R^3} \frac{\chi_i^{z_a}(\br-\bR_a)\, \chi_j^{z_b}(\br-\bR_b) \, \chi_k^{z_c}(\br'-\bR_c)\, \chi_\ell^{z_d}(\br'-\bR_d)}{|\br-\br'|} \, \d\br \, \d\br',
\end{align*}
arising in discretized electronic structure models can then be computed
analytically by means of explicit calculations and recursion formulas.

However, individual Gaussian function poorly describes the cusps of the bound states electronic wave-functions at nuclear positions. \emph{Contracted} Gaussians \cite{mcweenyGaussianApproximationsWave1950}, that are linear combinations of Gaussians with different variances,
were quickly introduced
as they overcome this deficiency. Several classes of GTO basis sets have been proposed since the 50's: STO-$n$g basis sets \cite{Hehre1969-qn} were built as the contraction of $n$ Gaussians that fit Clementi STO SCF AOs in an $L^2$ least-squares sense~\cite{Stewart1969-hd}. It was quickly realized that better GTO basis sets could be obtained by minimizing atomic Hartree--Fock ground state energy. This approach led to the split-valence basis sets (e.g. 6-31G) with core and valence orbitals being approximated differently, developed by Pople et al. \cite{Binkley1980-bc}. Basis set better suited for correlated electronic structure methods were then introduced, notably Atomic Natural Orbitals (ANO) \cite{almlofGeneralContractionGaussian1987} and Dunning basis sets \cite{Dunning1989-sw}. ANO basis sets are built by selecting occupied and virtual orbitals from Hartree--Fock and natural orbitals from correlated computations of atomic systems.
Dunning bases provide a (finite) hierarchy of bases obtained by consistently increasing the number of basis functions corresponding to different angular momenta. This optimization strategy yields the so-called correlation consistent cc-pVXZ basis sets, which are, with their \emph{augmented} version, still commonly used nowadays.

\medskip

Mathematical studies proving convergence rates or proposing systematic
enrichment of GTO basis sets are so far quite limited. The approximability of
solutions to electronic structure problems by Gaussian functions was studied
in~\cite{Kutzelnigg1994-ys}, and later on in~\cite{Scholz2017-zv,Shaw2020-fh}.
An {\it a priori} error estimate on the approximation of Slater-type functions by Hermite and even-tempered Gaussian functions
was derived in~\cite{Bachmayr2014-zg}. A
construction of Gaussian bases combined with wavelets was proposed on a
one-dimensional toy model in~\cite{Pham2017-bd}.

\medskip

Commonly used Pople and Dunning GTO basis sets were optimized from atomic configuration
energies and Hartree--Fock (and/or natural) atomic orbitals.
In this article, we propose a different approach, which is adaptable to any criterion one might be
interested in, and involves molecular configurations. In \autoref{sec:math_setting}, we introduce an abstract
mathematical framework for the construction of optimal AO basis sets, based on the choices of
\begin{enumerate}
  \item a set of admissible atomic configurations $\Omega$;
  \item a probability measure $\PP$ on $\Omega$;
  \item a set of admissible AO basis sets $\mathcal{B}$;
  \item a criterion $j(\chi,\omega)$ quantifying the error between the exact values of the quantities of interest when the system has atomic configuration $\omega \in \Omega$ -- for the continuous model under consideration -- and the ones obtained after discretization in the basis set $\chi \in \mathcal{B}$.
\end{enumerate}
We also provide examples of possible choices of $\Omega$, $\PP$, $\mathcal{B}$,
and $j$. As a proof of concept (\autoref{sec:toy_model}), we apply this
strategy to a simple toy model of a 1D homonuclear diatomic \enquote{molecule}
with two 1D non-interacting spinless \enquote{electrons}, which allows for extremely
accurate reference calculations. Finally, we present numerical results in
\autoref{sec:num}, where we compare the efficiency of the so-optimized AO
bases compared to AO basis constructed from the occupied and unoccupied orbitals of the isolated \enquote{atom}.

\section{Optimization criteria} \label{sec:math_setting}

\subsection{Abstract framework}

We start by formulating the problem of basis set optimization in an abstract setting. The procedure can be divided into four steps.

\medskip

First, we select the set $\Omega$ of all possible atomic configurations we are
{\it a priori} interested in. For instance, depending on the foreseen
applications, one can consider the set of all possible finite atomic
configurations containing only hydrogen, nitrogen, carbon, and oxygen atoms,
or the set of all possible periodic arrangements of chemical elements with
less than 20 atoms per unit cell.

\medskip

Second, we equip $\Omega$ with a probability measure $\PP$ in order to allow
for different configurations to have different weights in the optimization
procedure. We will see later that our method requires the computation of very
accurate reference solutions for all $\omega$'s in the support of $\PP$. For
practical reasons we therefore need to choose $\PP$ of the form
\begin{equation}\label{eq:generalP}
  \PP = \sum_{n =1}^{N_{\rm c}} w_n \delta_{\omega_n},
\end{equation}
where $\{\omega_1,\dots , \omega_{N_{\rm c}} \}$ is a finite (not too large) subset of $\Omega$, $\delta_{\omega_n}$ the Dirac mass at $\omega_n$, and $\{w_1,\dots , w_{N_{\rm c}} \}$ are positive weights such that $\sum_{n=1}^{N_{\rm c}} w_n=1$.  Assume that we are solely interested in reproducing accurately the dissociation curve of the HF (Hydrogen Fluoride) diatomic molecule. Then the set $\Omega$ should be identified with the interval $(0,+\infty)$, and a configuration $\omega \in \Omega$ with the H$-$F interatomic distance $R \in (0,+\infty)$, and $\PP$ should be a probability measure on the interval $(0,+\infty)$. The selection of the $\omega_n$'s and $w_n$'s can be done in various way. An option is to
\begin{enumerate}
  \item[i)]  choose a continuous probability distribution $\PP$ on $(0,+\infty)$ on the basis of chemical arguments, putting little weight on usually unimportant very small interatomic distances, more weight on  interatomic distances close to the equilibrium distance ($d \simeq 0.92$ \AA), sufficient cumulated weight on very large interatomic distance to correctly reproduce the dissociation energy, and more or less weight on intermediate interatomic distances in the range $2-8$ \AA, depending on its importance for the targeted application;
  \item[ii)] fix the number $N_{\rm c}$ according to the available computational means;
  \item[iii)] compute the $\omega_n$'s and $w_n$'s
    using e.g. quantization algorithms~\cite{merigot2021non} possibly based on optimal transport or clustering algorithms~\cite{Pages2015-vu}.
\end{enumerate}

\medskip

Third, we select the set $\mathcal{B}$ of admissible AO basis sets. Restricting ourselves to the framework of GTOs, this can be done by choosing, for each chemical element arising in $\Omega$, the number, symmetries, and contraction patterns of the Gaussian polynomials of the AO associated with this particular element. In this case, $\mathcal{B}$ has the geometry of a convex polyedron of $\R^d$.

\medskip

Given an atomic configuration $\omega \in \Omega$ and an AO basis set $\chi \in \mathcal{B}$, we denote by $\chi_\omega$ the one-electron finite-dimensional space obtained by using the AO basis set $\chi$ to describe the electronic structure of a molecular system with atomic configuration $\omega$ and an arbitrary number $N$ of electrons.

\medskip

The fourth and final step consists in choosing a criterion $j(\chi,\omega)$ quantifying the quality of the results obtained when using the basis set $\chi \in \mathcal{B}$ to compute the electronic structure of a molecular system with atomic configuration $\omega$. The choice of the function
\[
  j:\mathcal{B}\times\Omega \to \mathbb{R}_+
\]
depends on the quantity of interest (QoI) to the user, and on the respective weights of these quantities in the case of multicriteria analyses. For instance, if one focuses on the ground-state energy of the electrically neutral molecular system, a natural criterion is
\begin{equation}
  \label{eq:def_jE}
  j_E(\chi,\omega) \coloneqq \left| E_\omega - E_\omega^{\chi} \right|^2,
\end{equation}
where $E_\omega$ is the exact ground-state energy of the neutral system with atomic configuration $\omega$ for the chosen continuous model (e.g. Hartree--Fock, MCSCF, Kohn--Sham B3LYP\dots) and $E_\omega^\chi$ the ground-state energy obtained with the model discretized in the AO basis set $\chi$. Note that the absolute value of the difference is squared to make $j_E$ differentiable. Another possible choice is to use a criterion based on the one-body reduced density matrices (1-RDM), for instance
\begin{equation}
  \label{eq:def_jA}
  j_A(\chi,\omega) \coloneqq - {\rm Tr}\left(\Pi_{\chi_\omega}^A \gamma_\omega \Pi_{\chi_\omega}^A A\right),
\end{equation}
where $A$ is a given self-adjoint, positive, definite operator on the one-particle state space $\mathcal H$ with form domain $Q(A)$, $\gamma_\omega$ the exact ground-state 1-RDM of the neutral system with atomic configuration $\omega$ for the chosen continuous model, and $\Pi_{\chi_\omega}^A : Q(A) \to {\mathcal X}_\omega \subset \mathcal H$ the orthogonal projector on ${\mathcal X}_\omega$ for the inner product $A$ on $Q(A)$. If $A=I_{\mathcal H}$, then the $Q(A)=\mathcal H$ and $\Pi_{\chi_\omega}^A$ is the orthogonal projector on ${\mathcal X}_\omega$ for the inner product of $\mathcal H$. If $A=(1-\Delta)$, then $Q(A)$ is the Sobolev space $H^1(\R^3)$, and $\Pi_{\chi_\omega}^A$ is the orthogonal projector on ${\mathcal X}_\omega$ for the $H^1$-inner product. Diagonalizing $\gamma_\omega$ as
$$
\gamma_\omega = \sum_j n_{\omega,j} |\psi_{\omega,j}\rangle \langle \psi_{\omega,j}|, \quad 0 \lq n_{\omega,j} \lq 1, \quad \langle \psi_{\omega,j}| \psi_{\omega,j'} \rangle=\delta_{jj'},
$$
where the $n_{\omega,j}$'s are the natural occupation numbers (NON) and $\psi_{\omega,j}$'s the natural orbitals (NO) for the chosen continuous model of the neutral system with atomic configuration $\omega$, it holds
$$
j_A(\chi,\omega) = - \sum_j n_{\omega,j} \| \Pi_{\chi_\omega}^A \psi_{\omega,j}\|^2_{Q(A)}.
$$
Minimizing $j_A(\chi,\omega)$ thus amounts to maximizing the NON-weighted sum of the $Q(A)$-norms of $Q(A)$-orthogonal projections of the NON on the discretization space ${\mathcal X}_\omega$. Other criteria may include errors on molecular properties, or a weighted sum of several elementary criteria, each of them targeting a specific property. The criterion should be chosen according to the intended application.

\medskip

The aggregated criterion to be optimized then reads as an integral over the configuration space $\Omega$ with respect to the probability measure $\PP$
\begin{equation}\label{eq:def_J}
  J(\chi) \coloneqq \int_\Omega j(\chi, \omega) \d\PP(\omega),
\end{equation}
and the problem of basis set optimization can be formulated as
\[
  \boxed{{\rm  find\;}\chi_0 \in \underset{\chi\in\mathcal{B}}{\rm argmin\;} J(\chi)}
\]

\begin{remark}[Reference solutions]
  The evaluation of criteria $J_E$ and $J_A$ hinges on the knowledge of exact GS energy $E_\omega$ or 1-RDM $\gamma_\omega$ for $\omega$ in the support of $\PP$. In practice, these data can be approximated by very accurate off-line reference computations for a small, wisely chosen, sample of configurations $\omega$. This is the reason why the probability measure $\PP$ can only be a finite weighted sum of Dirac distributions, as defined in \eqref{eq:generalP}.
\end{remark}

\section{Application to 1D toy model} \label{sec:toy_model}

In this section, we focus on a linear one-dimensional toy model, mimicking a homonuclear diatomic molecule.

\subsection{Description of the model}

\begin{figure}[h!]
  \centering
  \includegraphics[scale=.8]{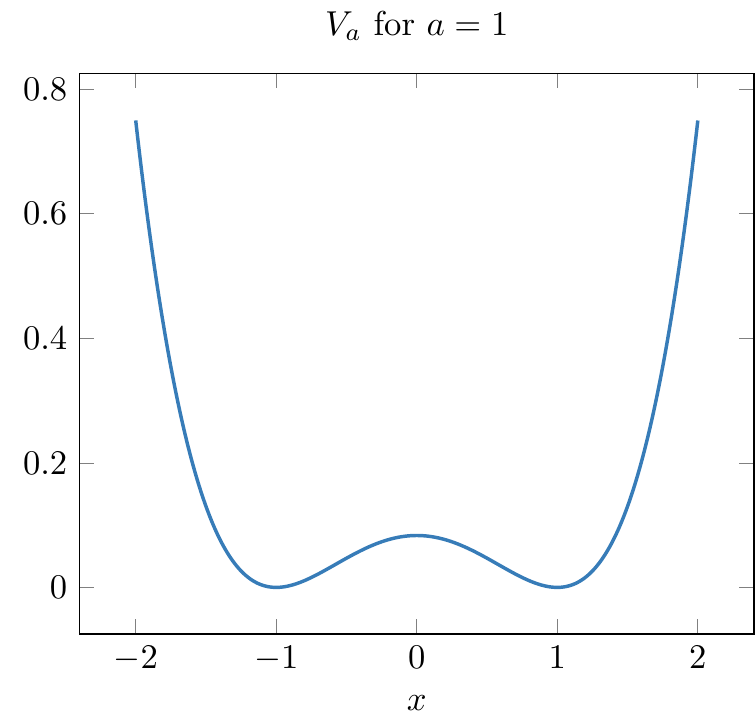}
  \phantom{aaaaaaa}
  \includegraphics[scale=.8]{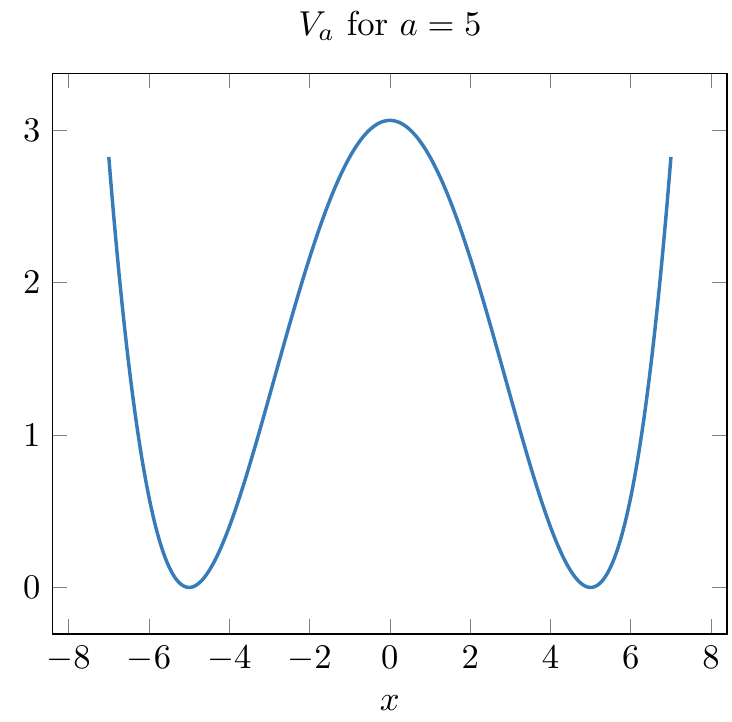}
  \caption{$x\mapsto V_a(x)$ for $a=1$ and $a=5$.}
  \label{fig:toy_model}
\end{figure}

Let us consider a system of two 1D point-like \enquote{nuclei} and two 1D spinless non-interacting quantum \enquote{electrons}. The one-particle state space is then ${\mathcal H}=L^2(\R)$ and the configuration space $\Omega =(0,+\infty)$. In this section, a configuration of $\Omega$ will be labelled by the positive real number $a >0$ such that the nuclei are located at $-a$ and $a$. The one-particle Hamiltonian at configuration $a$ then is
\begin{equation}\label{eq:toy_model_hamiltonian}
  H_a = -\frac{1}{2}\frac{{\rm d}^2}{{\rm d}x^2} + V_{a},
\end{equation}
where $V_a$ models the nuclei-electron interaction. We choose $V_{a}$ to be a double-well potential with minima at $-a$ and $+a$, defined by
\begin{equation}
  \Forall x\in\mathbb{R},\quad V_{a}(x)=\frac{1}{8a^2 + 4}(x-a)^2(x+a)^2.
\end{equation}
Several considerations led us to define the potential as such. First, $V_a$ is designed so that i) each $H_a$ admits a non-degenerate ground-state of energy $E_a$, and ii) the function $a \mapsto E_a$ has the shape of the ground-state dissociation curve of a homonuclear diatomic molecule with atoms at $-a$ and $+a$. Since the two \enquote{electrons} do not interact, the ground-state energy $E_a$ and density matrices $\gamma_a \in \mathcal{G}_2$ are given by
\begin{equation}\label{eq:pb_init}
  E_a = {\rm Tr}\left(H_a\gamma_a \right) = \underset{\gamma\in\mathcal{G}_2}{\rm min}\,{\rm Tr}\left(H_{a}\gamma\right) ,
\end{equation}
where
\begin{equation*}
  \mathcal{G}_2 \coloneqq \left\{ \gamma \in \mathcal{L}(L^2(\RR)),\ \gamma^2=\gamma=\gamma^*,\ {\rm Tr}(\gamma) = 2 \right\},
\end{equation*}
$\mathcal{L}(L^2(\RR))$ denoting the space of bounded linear operators on $L^2(\RR)$. The existence and uniqueness of the solution to problem~\eqref{eq:pb_init} can be shown by elementary arguments of functional analysis and spectral theory that we do not detail here. Second, $V_0(x)=\tfrac{1}{4}x^4$ so that \eqref{eq:toy_model_hamiltonian} corresponds to the quartic oscillator, for which we have reference numerical solutions (e.g. \cite{blinder2019eigenvalues}).
Third, $V_a$ behaves like $x^2/2$ around $\pm a$ for large values of $a$ and $V_a(0)\sim a^4/8\to+\infty$ when $a\to+\infty$. Therefore, in the limit $a\to +\infty$,  problem \eqref{eq:pb_init} is similar to two decoupled quantum harmonic oscillators centered in $-a$ and $+a$ whose bound states are all explicitly known. For the sake of illustration, we display in \autoref{fig:toy_model} the potential $V_a$ for two different values of $a$.

\medskip
In practice, it is convenient to compute $\gamma_a$ and $E_a$ from the lowest two eigenvalues $\lambda_{a,1} < \lambda_{a,2}$ of $H_a$ and an associated pair $(\varphi_{a,1},\varphi_{a,2})$
of orthonormal eigenvectors
\begin{equation}
  \begin{cases}
    \displaystyle
    H_a \varphi_{a,i} = \lambda_{a,i} \varphi_{a,i}, \quad i=1,2 \label{eq:pb_init_2} \\
    \displaystyle
    \langle \varphi_{a,i}| \varphi_{a,j}\rangle = \delta_{ij}, \quad i,j=1,2,
  \end{cases}
\end{equation}
$\langle \cdot | \cdot \rangle$ denoting the $L^2$ inner product. We indeed have
\begin{equation}
  E_a = \lambda_{a,1} + \lambda_{a,2}\quad\mbox{and}\quad \gamma_a=|\varphi_{a,1}\rangle\langle\varphi_{a,1}|+|\varphi_{a,2}\rangle\langle\varphi_{a,2}|.
\end{equation}

The evaluation of criteria $J_A$ and $J_E$ requires the computation of reference ground-state density matrices or energies, which amounts to find very accurate solutions of \eqref{eq:pb_init_2} for the configurations $a_k$ in the support of the chosen atomic probability measure
\begin{equation}\label{eq:propP}
  \PP = \sum_{n=1}^{N_{\rm c}} w_n \delta_{a_n}, \quad 0 < a_1 < a_2 < \cdots
  < a_{N_{\rm c}}, \quad w_n > 0, \quad \sum_{n=1}^{N_{\rm c}} w_n = 1.
\end{equation}
We chose to compute these reference data using a 3-point
finite-difference (FD) scheme on a large enough interval $[-x_{\rm max}, x_{\rm max}]$ discretized into a uniform grid with $N_g$ grid points:
$$
x_j=-x_{\rm max}+j \delta x, \quad 1 \lq j \lq N_g, \quad \delta x=\frac{2 x_{\rm max}}{N_g+1}.
$$
We then impose homogeneous Dirichlet boundary conditions at $-x_{\rm max}$ and $x_{\rm max}$.
The parameter $x_{\rm max}$ is chosen such that $x_{\rm max}=a_{\rm max}+ r_{\rm max}$, where $a_{\rm max}=\max (\mbox{supp}(\PP))$ and $r_{\rm max} > 0$ is is the radius beyond which atomic densities are zero at machine (double) precision. Note that this numerical scheme is independent of the configuration $a$. The FD discretization of problem \eqref{eq:pb_init_galerkin} gives rise to the eigenvalue problem
\begin{equation}
  \begin{cases}
    \displaystyle
    H_{a}^{\rm FD}\varphi_{a,i}^{FD} = \lambda_{a,i}^{\rm FD}\varphi_{a,i}^{\rm FD} \quad i=1,2 \label{eq:pb_init_galerkin_FD} \\
    \delta x\bigl(\varphi_{a,i}^{\rm FD}\bigr)^T\varphi_{a,j}^{\rm FD}= \delta_{ij}^{\;},
  \end{cases}
\end{equation}
where $H_{a}^{\rm FD} \in \R^{N_g \times N_g}_{\rm sym}$ is a real symmetric matrix of size $N_g \times N_g$, and the reference data are obtained as
\begin{equation}\label{eq:discrete_GS_E_DM_FD}
  E_a^{\rm FD}=\lambda_{a,1}^{\rm FD}+\lambda_{a,2}^{\rm FD}\quad\mbox{and}\quad P_{a}^{\rm FD} = \varphi_{a,1}^{\rm FD}\bigl(\varphi_{a,1}^{\rm FD}\bigr)^T+\varphi_{a,2}^{\rm FD}\bigl(\varphi_{a,2}^{\rm FD}\bigr)^T\in\mathbb{R}^{N_g\times N_g}_{\rm sym},
\end{equation}
where $P_{a}^{\rm FD}$ can be interpreted as an approximation of the matrix $\gamma_{a}(x_j,x_{j'})$ containing the values of the (integral kernel of the) density matrix $\gamma_a$ at the grid points.

\subsection{Variational approximation in AO basis sets}\label{subsec:discrete_setting}

For any given configuration $a\in\mathbb{R}_+$ and basis $\chi=\{\chi_\mu\}_{1\lq \mu\lq N_b}\in\mathcal{B}$,  problem \eqref{eq:pb_init_2} is solved using a Galerkin method with the basis $\chi_a=\{\chi_{a,\mu}\}_{1 \lq \mu \lq 2N_b}$ composed of two copies of the basis $\chi$, the first one translated to $a$, and the second one to $-a$:
$$
\chi_{a,1}=\chi_1(\cdot-a), \dots, \chi_{a,N_b}=\chi_{N_b}(\cdot-a), \; \chi_{a,N_b+1}=\chi_1(\cdot+a), \dots, \chi_{a,2N_b}=\chi_{N_b}(\cdot+a).
$$
Defining the Hamiltonian matrix
\[
  H_{a}^{\chi} = \left(\left\langle\chi_{a,\mu} \biggl|\left(-\frac{1}{2}\frac{\d^2}{\d x^2}+V_{a}\right) \biggr| \chi_{a,\nu}\right\rangle\right)_{1\lq\mu,\nu\lq 2N_b}
\]
and the overlap matrix
\[
  S^{\chi}_a=\left(\langle\chi_{a,\mu}|\chi_{a,\nu}\rangle\right)_{1\lq\mu,\nu\lq 2N_b},
\]
the discretization of problem \eqref{eq:pb_init_2} in the AO basis set $\chi$ then reads as the generalized eigenvalue problem: find $\bigl(C_{a,i}^{\chi},\lambda_{a,i}^\chi\bigr)\in\mathbb{R}^{2N_b}\times\mathbb{R}$, $i=1,2$ such that
\begin{equation}
  \label{eq:pb_init_galerkin}
  \begin{cases}
    \displaystyle
    H_{a}^{\chi}C_{a,i}^{\chi} = \lambda_{a,i}^{\chi}S^{\chi}_a C_{a,i}^{\chi} \quad i=1,2 \\
    \bigl(C_{a,i}^{\chi}\bigr)^TS^{\chi}_a C_{a,j}^{\chi}= \delta_{ij}.
  \end{cases}
\end{equation}
The approximation $\varphi_{a,i}^{\chi}$ of $\varphi_{a,i}^{\,}$ in the AO basis set $\chi$ can then be recovered as the linear combination of atomic orbitals (LCAO)
\begin{equation}
  \forall x\in\mathbb{R},\quad \varphi_{a,i}^{\chi}(x) = \sum\limits_{\mu=1}^{2N_b}[C_{a,i}^{\chi}]_\mu^{\,}\chi_{a,\mu}^{\,}(x).
\end{equation}
One way to compare the LCAO ground-state 1-RDM to the reference FD solution
$P_a^{\rm FD}$ is to simply evaluate the former at the grid points $x_j$,
which gives rise to the matrix $P_a^\chi  \in \R^{N_g \times N_g}_{\rm sym}$ with entries
$$
[P_a^\chi]_{jj'} = \sum_{i=1}^2 \varphi_{a,i}^{\chi}(x_j) \varphi_{a,i}^{\chi}(x_{j'}).
$$
Due to numerical errors, the matrix $P_a^\chi$ is however not a rank-2 orthogonal projector. We therefore chose to follow a slightly different route (leading to very similar results). The finite difference grid gives a reference discrete setting in which any quantity of interest for any configuration and AO basis set can be expressed. For all $a$'s, the basis $\chi_a$ is represented by a matrix $X_a\in\mathbb{R}^{N_g\times 2N_b}$.
For any vectors $Y_1,Y_2\in\mathbb{R}^{N_g}$, the discrete $A$ inner product simply reads $\delta x Y_1^T A Y_2^{\,}$ where the notation $A$ stands for both the continuous inner product and its finite-difference discretization matrix. We denote by $\Vert\cdot\Vert_A$ the associated norm on $\R^{N_g}$.
Solutions to \eqref{eq:pb_init_galerkin} are then obtained by
approximating respectively the Hamiltonian and overlap matrix by
\[
  H_a^{\chi} \simeq H_a^{X} \coloneqq
  \left( \delta x X_{a,\mu}^T H_a^{\rm FD} X_{a,\nu}^{\,}
  \right)_{1\lq\mu,\nu\lq 2N_b}, \quad
  S^{\chi}_a \simeq S^X_a \coloneqq
  \left( \delta x X_{a,\mu}^T X_{a,\nu}^{\,}
  \right)_{1\lq\mu,\nu\lq 2N_b},
\]
and finding $( C_{a,i}^{X},\lambda_{a,i}^{X})\in \RR^{2N_b}\times\RR$, $i=1,2$, such that
\begin{equation}
  \left\{
    \begin{array}{l}
      \displaystyle
      H_a^{X}  C_{a,i}^{X} = \lambda^{X}_{a,i} S^{X}_a  C_{a,i}^{X},\quad i=1,2 \label{eq:discrete_galerkin} \\ \\
      (C_{a,i}^{X})^T S^{X}_a C_{a,j}^{X} = \delta_{ij}^{\;},\quad i,j=1,2,
    \end{array}
  \right.
\end{equation}
from which we get the discrete approximations
\begin{equation}
  \label{eq:discrete_phi}
  {\varphi}_{a,i}^X = X_a^{\,} C_{a,i}^X,\quad i=1,2.
\end{equation}
Let us gather the coefficients ${C}_{a,i}^X$ into the $2N_b\times 2$ matrix
${C}_a^X=\left(\,C_{a,1}^X\,|\,C_{a,2}^X\,\right)$.
The ground-state density matrix in the basis $\chi_a$ is approximated by
\begin{equation}
  \label{eq:discrete_GS_DM}
  {P}_a^{X}={\varphi}_{a,1}^{X}({\varphi}_{a,1}^{X})^T+{\varphi}_{a,2}^{X}({\varphi}_{a,2}^X)^T = \left(X_a^{\,} {C}^X_a\right) \left(X_a^{\,} {C}^X_a\right)^T \in \RR^{N_g\times N_g}_{\rm sym}.
\end{equation}

\subsection{Overcompleteness of Hermite Basis Sets}
\label{subsec:HBS}

\begin{figure}[!h]
  \centering
  \includegraphics[width=0.5\linewidth]{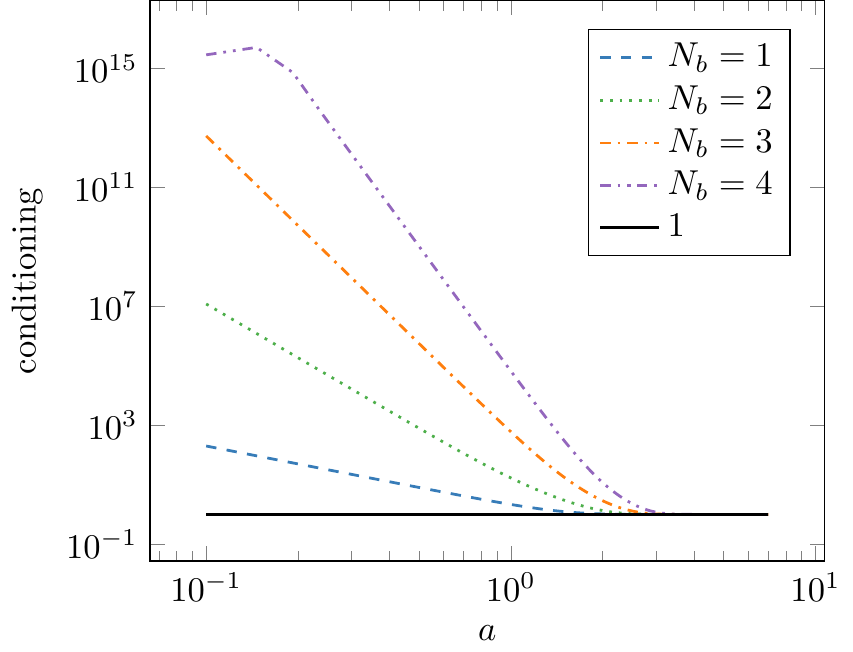}
  \caption{Condition number of the HBS overlap matrix $S_a^{\chi^{\rm HBS}} $ for different values of $a$ in log-log scale. The larger the basis set, the faster the condition number blows up for small values of $a$.}
  \label{fig:cond_S}
\end{figure}

Before getting into basis set optimization, we introduce the following standard Hermite Basis Set (HBS), constructed from eigenfunctions of the quantum harmonic oscillator. Those functions are solutions to the eigenvalue problem $\left(-\frac{1}{2}\frac{{\rm d}^2}{{\rm d}x^2} + \frac12 x^2 \right)h_n = \varepsilon_n h_n$ and are explicitly given by
\begin{equation}\label{def:HBS}
  h_n(x) = c_np_n(x)\exp\left({-\frac{x^2}{2}}\right),\quad \varepsilon_n=n+\frac{1}{2}, \quad n \in \NN,
\end{equation}
where $p_n$ is the Hermite polynomial of degree $n$ (with the same parity as $n$) and $c_n$ a
normalization constant such that $(h_n)_{n\in\NN}$ forms an orthonormal basis of $L^2(\RR)$.
The $h_n$'s are the analogues of the standard atomic orbitals obtained by solving atomic electronic structure problems. Let us first consider the AO basis set made of the first $N_b$ Hermite functions
$$
\chi^{\rm HBS}= \{ \chi_\mu^{\rm HBS} \}_{1 \lq \mu \lq N_b} = \{h_n\}_{0 \lq n \lq N_b-1}.
$$
The overlap matrix for the configuration $a$ then is of the form
\begin{equation*} S_{a}^{\chi^{\rm HBS}} =
  \begin{pmatrix} I_{N_b} & \Sigma_{a} \\[0.1cm] \Sigma_{a}^T & I_{N_b} \end{pmatrix}
  \quad\text{where}\quad
  \Sigma_{a} \coloneqq (\langle h_{n}(\cdot-a)| h_{m}(\cdot+a)\rangle)_{0\lq n,m \lq N_b-1}.
\end{equation*}
The matrix $\Sigma_{a}$ corresponds to the overlap of functions that are
localized at different atomic positions. It satisfies $\Sigma_{a} \simeq 0$
when $a$ is large and $\Sigma_{a} \simeq I_{N_b}$ when $a$ is close to 0,
therefore causing conditioning issues on the overlap matrix
$S_{a}^{\chi^{\rm HBS}}$, a phenomenon known as \emph{overcompleteness}: when
$a$ is too small, the basis functions centered at~$\pm a$ are almost equal,
hence almost linearly dependent in the basis set. We illustrate this problem
by plotting the condition number of the overlap matrix $S_a^{\chi^{\rm HBS}}$ for different
values of $a$ in \autoref{fig:cond_S}, which indeed blows up for small values
of $a$.

\subsection{\texorpdfstring{Practical computation of the criterion $J_A$ and $J_E$}{}}

The rest of this section is dedicated to the rewriting and the computation of criteria $J_A$ and $J_E$ for our 1D model in the discrete setting.

\subsubsection{Reference orthonormal basis}

In order to avoid potential numerical stability issues, each of the $N_b$ atomic orbital $\chi_\mu$ is decomposed on a given truncated orthonormal basis of~$L^2(\RR)$ of size $\Nc$ such that $N_b \ll \Nc \ll N_g$. We choose here the orthonormal basis introduced in \eqref{def:HBS}.
Hence, the matrix $X_a \in \RR^{N_g \times 2N_b}$ is written as
\begin{equation}
  \label{def:X_a}
  X_a = B_aI_R,
\end{equation}
with
\begin{equation*}
  B_a = \begin{pmatrix}
    h_0(x_{\cdot}-a) | \cdots | h_{\Nc-1}(x_{\cdot}-a) | h_0(x_{\cdot}+a) | \cdots | h_{\Nc-1}(x_{\cdot}+a)
  \end{pmatrix}
  \in\mathbb{R}^{N_g\times 2\Nc},
\end{equation*}
and
\begin{equation}
  I_R = \begin{pmatrix}
    R & 0 \\ 0 & R
  \end{pmatrix}
  \in\mathbb{R}^{(2\Nc)\times (2N_b)},
\end{equation}
where $R\in\RR^{\Nc\times N_b}$ gathers the coefficients of the atomic orbitals $\chi_\mu$ in the truncated HBS orthonormal basis. Note that we have duplicated $R$ in $I_R$ as we consider the same basis at each position $\pm a$, but everything that follows can be easily adapted to the case where we would like to optimize the bases at each position separately (to deal with heteronuclear molecular systems for instance).
We moreover impose that  $R^T R = I_{N_b}$, so that the overlap matrix of $X_a$, denoted by $S(X_a)$, has the same form as in \autoref{subsec:HBS}, that is
\begin{equation}
  S(X_a) \coloneqq \delta x X_a^TX_a^{\,} = \begin{pmatrix} I_{N_b} & \Sigma_a \\ \Sigma_a^T & I_{N_b} \end{pmatrix},
\end{equation}
where $\Sigma_a$ is the overlap between functions localized at $+a$ and functions localized at $-a$. To avoid any issues arising from the conditioning of $S({X_a})$, the minimal sampled distance $a_{\rm min}$ should not be taken too small.

\medskip

In the following, we detail the computation of each of the two criteria using the matrix $R$ as the main variable. We will subsequently optimize the criteria $J_A$ and $J_E$ with respect to $R$ to obtain optimal AO basis sets. In order to ease the reading of the following computations, every vector of $\mathbb{R}^{N_g}$ is rescaled by a factor $\sqrt{\delta x}$ so that for any given $Y_1, Y_2\in\mathbb{R}^{N_g}$ the discrete $A$ inner product simply reads $Y_1^TAY_2^{\,}$. The same holds for overlap matrices: with this convention, $S(X_a) = X_a^TX_a^{\,}$. The output of the optimization is then scaled back to its former state by a factor $1/\sqrt{\delta x}$ to recover the original normalization.

\subsubsection{\texorpdfstring{Criterion $J_A$}{}}

Let $a\in\mathbb{R}_+$ be fixed and let $S^A(Y) = Y^T A Y$ denote the overlap matrix for the $A$-inner product of any rectangular matrix~$Y \in \R^{N_g \times d}$.
Since the columns of ${X_a}[S^A(X_a)]^{-\frac{1}{2}}$
are orthonormal for the $A$ inner product, that is
\[
  \left({X_a}[S^A(X_a)]^{-\frac{1}{2}}\right)^T
  A
  \left({X_a}[S^A(X_a)]^{-\frac{1}{2}}\right) = I,
\]
the projection $\Pi_{X_a}^{A}$ takes the simple form
\begin{equation}\label{eq:discrete_projection}
  \Pi_{X_a}^{A} = \left({X_a}[S^A({X_a})]^{-\frac{1}{2}}\right)\left({X_a}[S^A({X_a})]^{-\frac{1}{2}}\right)^T A = X_a [S^A({X_a})]^{-1}X_a^T A.
\end{equation}
Hence, using the cyclicity of the trace and definitions \eqref{eq:def_jA}, \eqref{eq:discrete_GS_E_DM_FD} and \eqref{eq:discrete_projection}, one has
\begin{align*}
  j_A(\chi,a) & \simeq
  - {\rm Tr}\left(P_a^{\rm FD}\, \Pi_{X_a}^{A}A \Pi_{X_a}^{A}\right)\\
  & = -{\rm Tr}\left(P_a^{\rm FD}\times (AB_aI_R^{\,})[S^A(B_aI_R^{\,})]^{-1}(AB_aI_R^{\,})^T\right)\\[0.15cm]
  & = -{\rm Tr}\left(M^{\rm offline}_A(a)I_R^{\,}[S^A(B_aI_R^{\,})]^{-1} I_R^T\right),
\end{align*}
where we have collected in the last expression all matrices independent of $R$ into the matrix
\begin{equation}
  M^{\rm offline}_{A}(a)=(AB_a)^TP_a^{\rm FD}AB_a \in \RR^{2\Nc\times 2\Nc}.
\end{equation}
Then, using the probability measure $\PP$ in \eqref{eq:propP}, we get
\[
  J_A(R) =  -\int_\Omega {\rm Tr}\left(M^{\rm
      offline}_{A}(a)I_R^{\,}[S^A(B_a^{\,}I_R^{\,})]^{-1} I_R^T\right) \,
  {\rm d}\PP(a) = - \sum_{n=1}^{N_{\rm c}} w_n {\rm Tr}\left(M^{\rm
      offline}_{A}(a_n)I_R^{\,}[S^A(B_{a_n}^{\,}I_R^{\,})]^{-1}
    I_R^T\right)
\]
and the optimization problem finally writes, with unknown $R\in \RR^{\Nc \times N_b}$ and for a given inner product $A$
\begin{equation}\label{pb:discrete_J_A}
  \boxed{\mbox{Find } R_{\rm opt}\in\underset{R\in\mathbb{R}^{\Nc\times N_b},\ R^TR=I_{N_b}}{\rm argmin} J_{A}(R)}
\end{equation}

\subsubsection{\texorpdfstring{Criterion $J_E$}{}}

Let again $a\in \mathbb{R}_+$ be fixed. We denote by
\[
  G(N_g,2) \coloneqq \{P\in\mathbb{R}^{N_g\times N_g}\,|\, P^2=P=P^T, {\rm Tr}(P)=2\}
\]
the discrete counterpart of the Grassmann manifold $\mathcal{G}_2$, and write $E_a^R$ (resp. $H_a^R$) instead of $E_a^{\chi}$ (resp. $H_a^\chi$), so that the dependence in the matrix $R$ appears explicitly. Equation \eqref{eq:pb_init} reads in the discrete setting
\begin{equation}
  \label{eq:E_crit_JE}
  \begin{split}
    E_{a}^R = \displaystyle{\min_{P\in \,G(N_g,2)}{\rm Tr}\left(P H_a^{R}\right)}
    & = \displaystyle{\min_{\substack{C\in\mathbb{R}^{2N_b\times 2}\\ (C)^TS(B_aI_R)C = I_2}} {\rm Tr}\left(CC^T\times (B_a I_R)^T H_a^{\rm FD} (B_a I_R)\right)} \\
    &={\rm Tr}\left(C_a^{R}(C_a^R)^T \times I_R^T M^{\rm offline}_E(a)I_R^{\,}\right)
  \end{split}
\end{equation}
where, as for the previous case, all matrices independent of $R$ have been gathered in the matrix
\begin{equation}
  M^{\rm offline}_E(a) = B_a^T H_a^{\rm FD} B_a^{\,},
\end{equation}
and the matrix $C_a^R$ is solution to the minimization problem
\begin{equation}\label{pb:discrete_E}
  \displaystyle{\min_{\substack{C^R\in\mathbb{R}^{2N_b\times 2}\\ (C^R)^TS(B_aI_R)C^R = I_2}} {\rm Tr}\left(C^R(C^R)^T\times I_R^T M_E^{\rm offline}(a)I_R \right)}
\end{equation}
and is given in practice by $C_a^R = \left[S(B_aI_R)\right]^{-\tfrac{1}{2}}(u_{a,1}|u_{a,2})$ where $u_{a,1}$ and $u_{a,2}$ are orthonormal eigenvectors associated to the lowest two eigenvalues of
\[
  \left[S(B_aI_R)\right]^{-\tfrac{1}{2}}I_R^{\,} M^{\rm offline}_{E}(a)I_{R}^T\left[S(B_aI_R)\right]^{-\tfrac{1}{2}}.
\]
From  \eqref{eq:propP} and \eqref{eq:E_crit_JE}, one can compute
\[
  J_E(R) = \int_{\Omega}\left|E_a^{\rm FD} - E_a^R\right|^2 \, {\rm d}\PP(a) =
  \sum_{n=1}^{N_{\rm c}} w_n \left|E_{a_n}^{\rm FD} - E_{a_n}^R\right|^2
\]
and the optimization problem reads
\begin{equation}\label{pb:discrete_J_E}
  \boxed{{\rm Find}\; R_{\rm opt} \in \underset{R\in\mathbb{R}^{\mathcal{N}\times N_b},\ R^TR=I_{N_b}}{\rm argmin}J_E(R)}
\end{equation}

\section{Numerical results}
\label{sec:num}

\subsection{Numerical setting and first results}
\label{subsec:num_res}

Problems \eqref{pb:discrete_J_A} and \eqref{pb:discrete_J_E} are solved by direct minimization algorithms over the Stiefel manifold~\cite{absil2009optimization}
\begin{equation*}
  {\rm St}(\Nc, N_b)=\{R\in\mathbb{R}^{\Nc\times N_b}\,|\, R^T R=I_{N_b}\}.
\end{equation*}
The explicit computation of the gradients of $J_A$ and $J_E$ with respect to
$R$ is detailed in the Appendix.  We used a L--BFGS algorithm (with tolerance
$10^{-7}$ on the norm of the projected gradient), as implemented in the {\it
  Optim.jl} package \cite{Optim.jl-2018} in the {\it Julia} language
\cite{Julia-2017}. As initial guess, we picked the first $N_b$ Hermite
functions introduced in~\autoref{subsec:HBS}.

\medskip

In this subsection, we choose a probability distribution $\PP$ supported in the interval $\Ic = [1.5,5]$ so as to retain the physics of interest that takes place around the equilibrium configuration $a_0\simeq 1.925$ and all the way to dissociation. In particular $a_{\rm min}= 1.5$ is taken sufficiently large to avoid the conditioning issues on the overlap matrices described in \autoref{subsec:HBS}. More precisely, all the results in this subsection are obtained with the probability
\begin{equation}\label{eq:probaP}
  \PP = \frac{1}{10} \sum_{n=1}^{10} \delta_{a_n} \quad \mbox{with } a_n=1.5+(n-1)\frac{3.5}9.
\end{equation}
The quantities $M^{\rm offline}_A(a_n)$ and $M^{\rm offline}_E(a_n)$ are computed offline beforehand. We will discuss this choice and consider other probability measures $\PP$ in Sections~\ref{subsubsec:extrapolation} and \ref{subsubsec:choiceP}.

The finite-difference grid is a uniform grid on the interval $[-20,20]$ discretized into $N_g = 1999$ points ($\delta x = 0.02$).
Finally, we decompose the basis functions to be optimized in the HBS $\{h_n\}_{0 \lq n \lq \Nc-1}$ of $L^2(\RR)$ of size $\Nc = 10$.
Regarding the choice of the inner product for the first criterion $J_A$, we used the standard $L^2(\RR)$ and the $H^1(\RR)$ inner products, and denoted $J_{L^2}$ and $J_{H^1}$ the corresponding. This translates at the discrete level by choosing $A=I_{N_g}$ for $J_{L^2}$ and $A = I_{N_g}-\Delta$ for $J_{H^1}$ where $\Delta$ is the 3-point finite-difference discretization matrix of the 1D Laplace operator.
Once obtained, the optimal bases are used to solve the variational problem \eqref{eq:discrete_galerkin} on a much finer sampling of $\Ic$ and their accuracy is compared to the HBS. The code performing the simulations and plotting the results is available online\footnote{\url{https://github.com/gkemlin/1D_basis_optimization}}. Also, for the sake of clarity in the plots, $\wt{E}_a$ (resp. $\wt{\rho}_a$) denotes the GS energy (resp. the density) in the configuration $a$ with a given basis (specified by the context) and $E_a$ (resp. $\rho_a$) stands for the reference energy (resp. density) on the finite difference grid. Note that we write HBS for the (nonoptimized) Hermite basis set, and $L^2$-OBS, $H^1$-OBS or $E$-OBS for optimized basis sets with respect to the criterion $J_{L^2}$, $J_{H^1}$, or $J_E$.

\medskip

\begin{figure}[h!]
  \centering
  \includegraphics[height=0.55\linewidth]{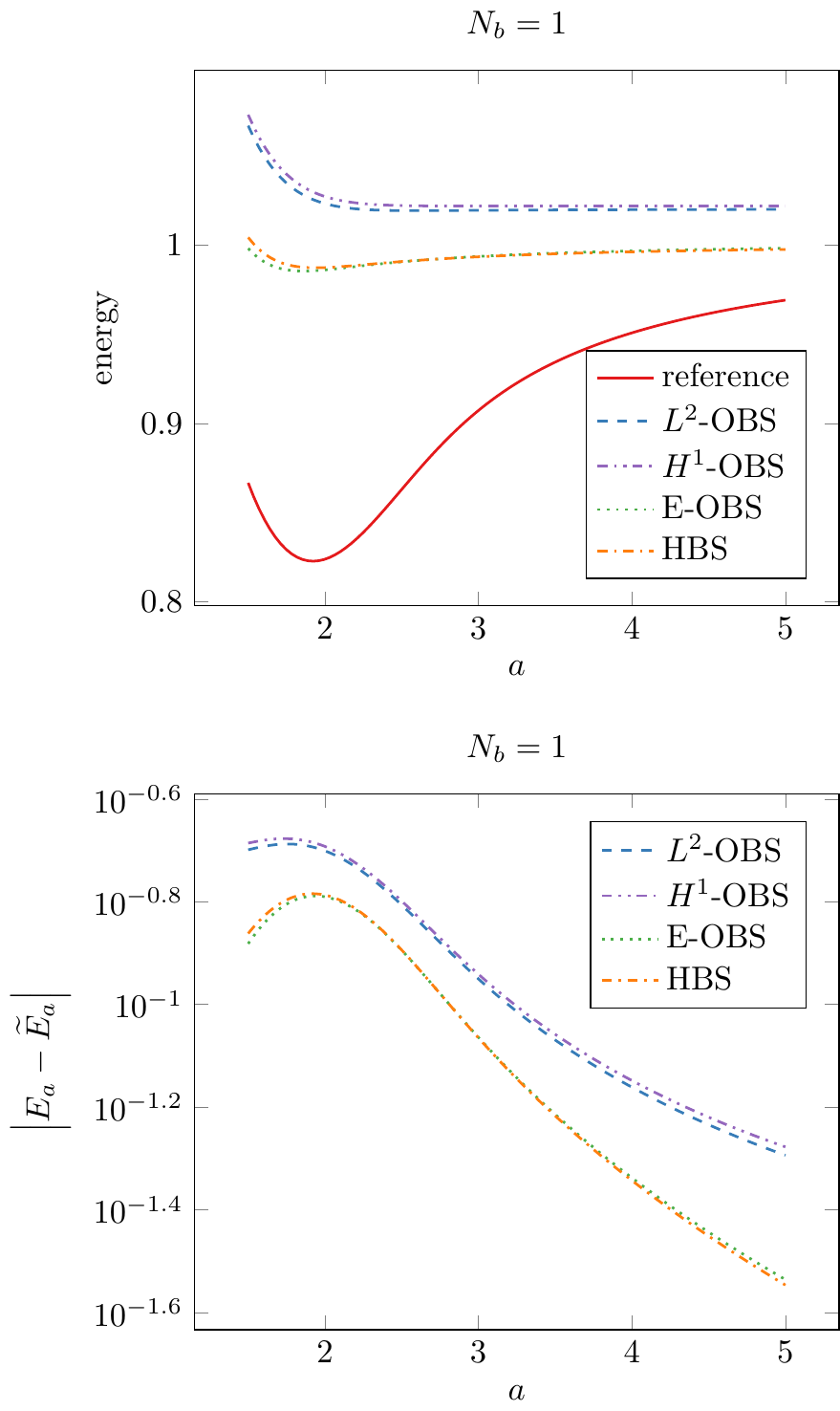}\hfill
  \includegraphics[height=0.55\linewidth]{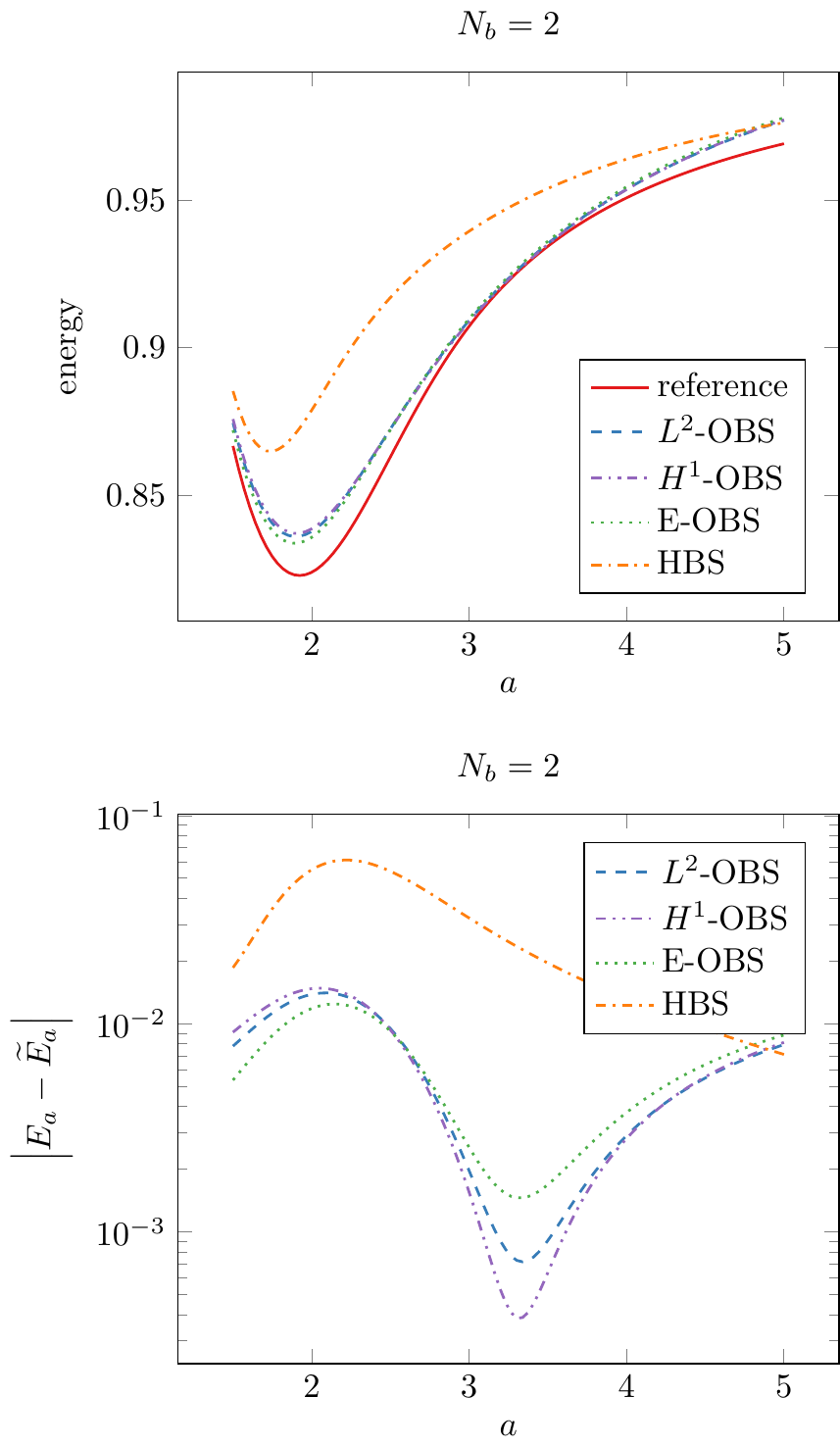}\hfill
  \includegraphics[height=0.55\linewidth]{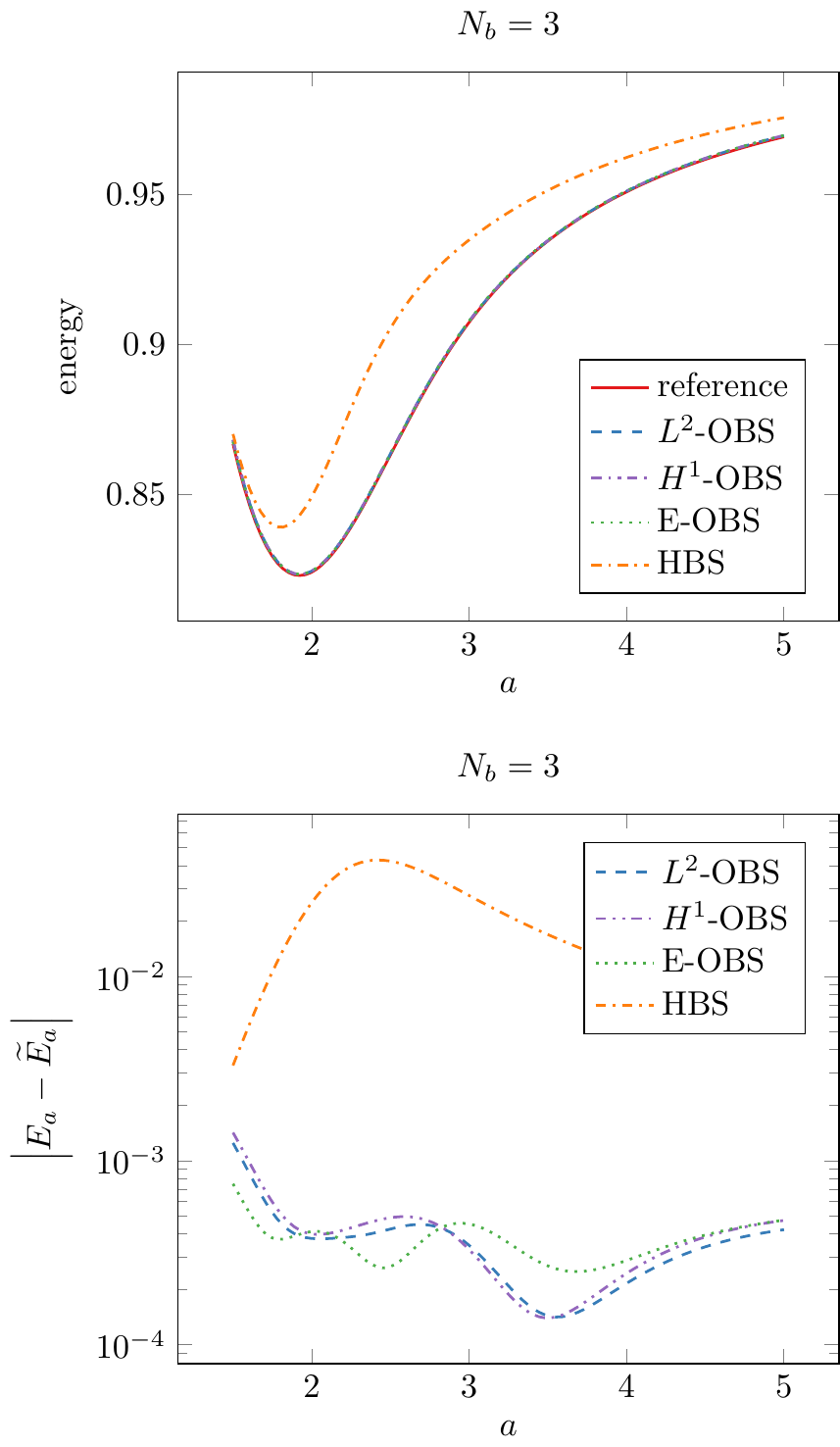}
  \caption{Energies for the optimal bases obtained with the different criteria. (Top) Dissociation curve. (Bottom) Errors on the energy on the range of configuration $\Ic=[1.5,5]$.}
  \label{fig:energy}
\end{figure}

\autoref{fig:energy} displays the dissociation curve and the energy difference on the interval $\Ic$ for different values of $N_b$, the size of the AO basis set.
For $N_b=1$, \ie\ only one basis function at $\pm a$, criterion $J_E$ shows better performance than the criterion $J_A$, regardless of the choice of norm to perform the projections. It also very closely matches the accuracy of the standard HBS.
When $N_b$ becomes larger however, the different criteria behave in a similar fashion and we observe that they approach the dissociation curve better than the Hermite basis. Comparing the values of criterion $J_E$ for all bases, which directly measures the distance to the dissociation curve, we see in \autoref{tab:criteria} that all optimized bases give an increased accuracy of roughly four orders of magnitude over the interval $\Ic$ for $N_b=4$.

\medskip

\begin{figure}[h!]
  \centering
  \hfill
  \includegraphics[width=0.3\linewidth]{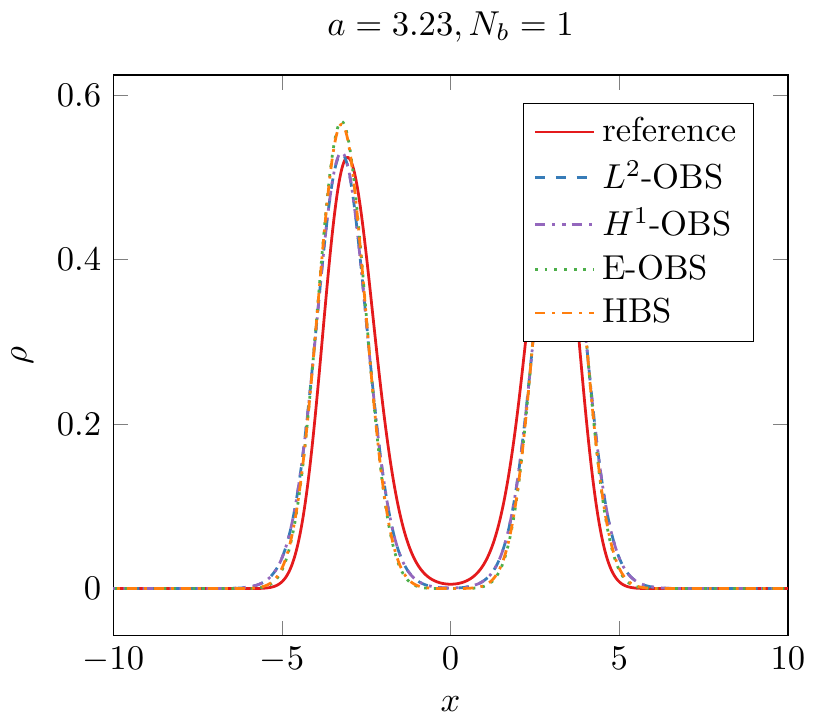}\hfill
  \includegraphics[width=0.3\linewidth]{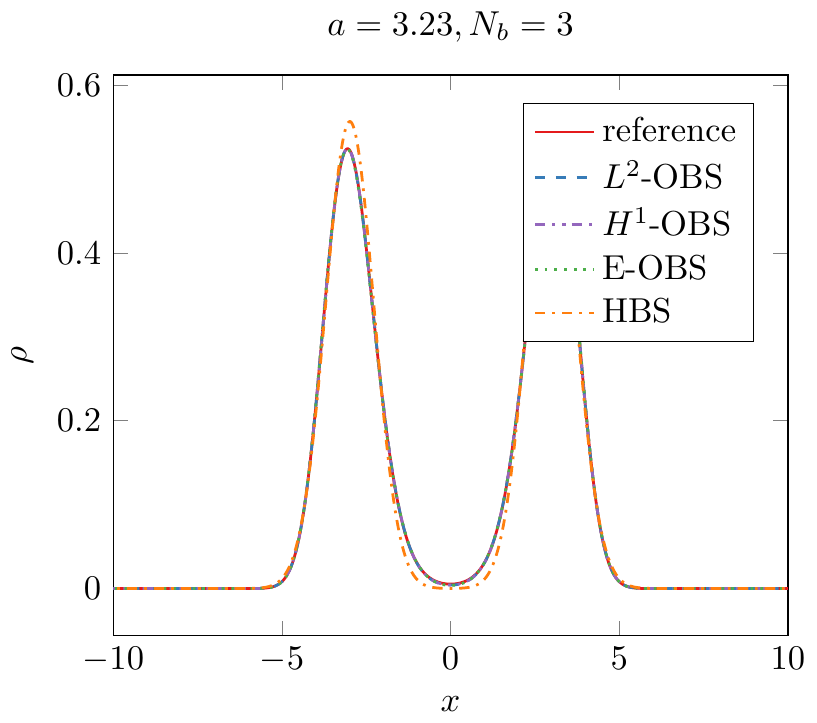}\hfill{\phantom{e}}\\
  \includegraphics[width=\linewidth]{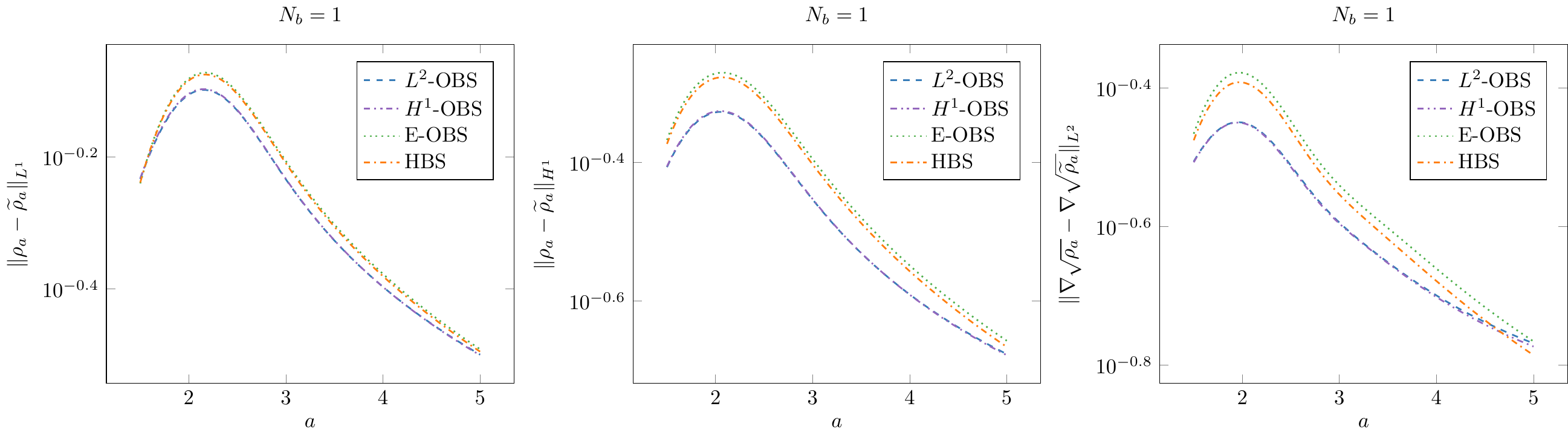}
  \includegraphics[width=\linewidth]{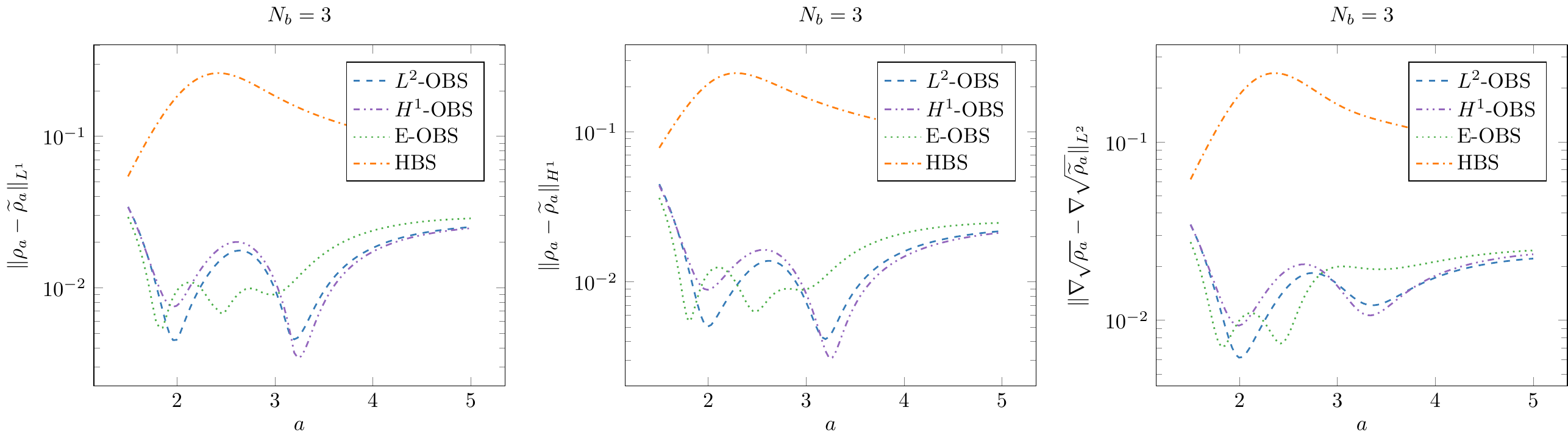}
  \caption{(Top) Densities for the optimal bases obtained with the different criteria. (Middle) Errors on the density for different norms with $N_b = 1$. (Bottom) Error on the density for different norms with $N_b = 3$.
  }
  \label{fig:density}
\end{figure}

In \autoref{fig:density}, we plot the density for a given value of $a$ and the error on the density for different norms, with varying values of $N_b$. The error is plotted with respect to three different distances: the $L^1$-norm, which corresponds to the $L^2$-norm on eigenvectors, the $H^1$-norm of the error on the density, as it is common to compute the forces $\int_\RR \rho \nabla_a V_a$ with good estimates on the $H^{-1}$-norm of $\nabla_aV_a$ (see e.g. \cite{cancesPracticalErrorBounds2021}), and the distance
\begin{equation*}
  \Vert \nabla\sqrt{\rho_1} - \nabla\sqrt{\rho_2} \Vert_{L^2}
\end{equation*}
(recall that the von-Weizs\"acker kinetic energy reads $\frac 12 \int_\RR |\nabla\sqrt{\rho}|^2$).
We observe similar behaviors between these different distances.
For $N_b=1$, both bases obtained with the first criterion behave slightly better than the standard Hermite basis and the basis computed with the second criterion. For $N_b = 3$, we observe again that all optimal bases yield better accuracy than the Hermite basis.
\autoref{tab:criteria} gives the confirmation that each basis for a given criterion indeed performs better than the other bases for that particular criterion. As for dissociation curves, we read from the values of $J_{L^2}$ and $J_{H^1}$  that the optimized bases yield similar results for large $N_b$, all of them giving lower values than the HBS. Note that the optimal bases for criterion $J_{L^2}$ and $J_{H^1}$ give similar results for any number of basis functions $N_b$, so that the $L^2$ and $H^1$ norm optimizations seem equivalent.

\medskip

In terms of computational time, first note that criterion $J_{H^1}$ is always more expensive to compute than $J_{L^2}$ as it requires additional matrix-vector products with the matrix $A$, this having noticeable impact on the computational time. Second, criterion $J_E$ requires less off-line data as it only needs to be given the reference eigenvalues while criterion $J_A$ requires the reference GS eigenvectors (or density matrices). This also influences the computational time. Indeed, regarding timings, criterion $J_E$ is about 10 times faster to minimize than criterion $J_{L^2}$ norm in our implementation.

\medskip

\begin{figure}[p]
  \centering
  \begin{tabular}{c}
    $N_b = 1$ \\
    \includegraphics[width=0.25\linewidth]{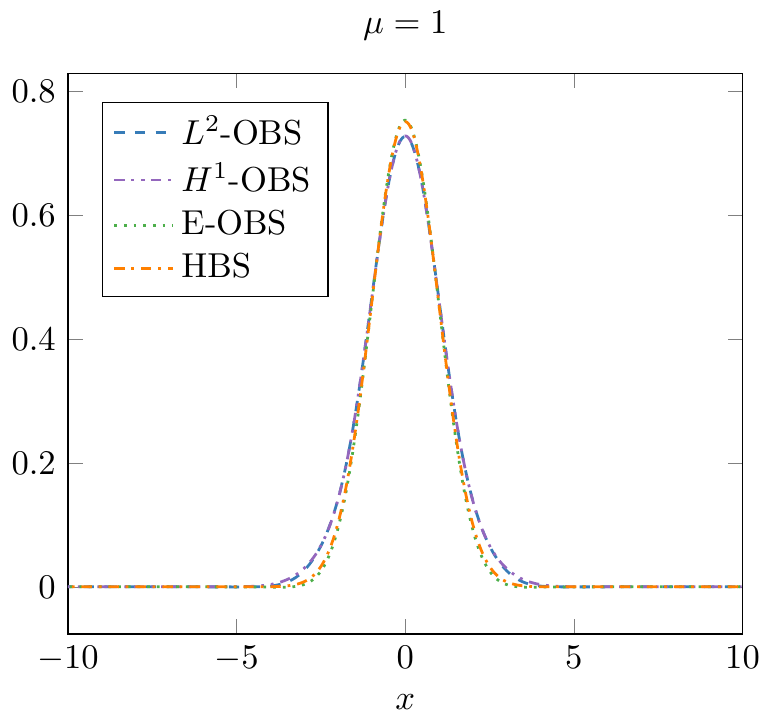}
    \\ \\
    $N_b = 2$ \\
    \includegraphics[width=0.25\linewidth]{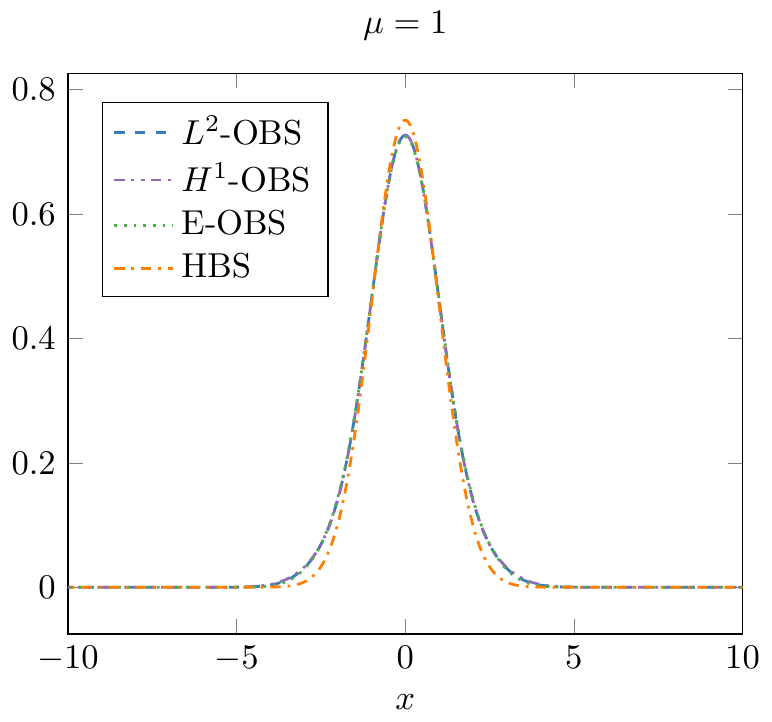}\includegraphics[width=0.25\linewidth]{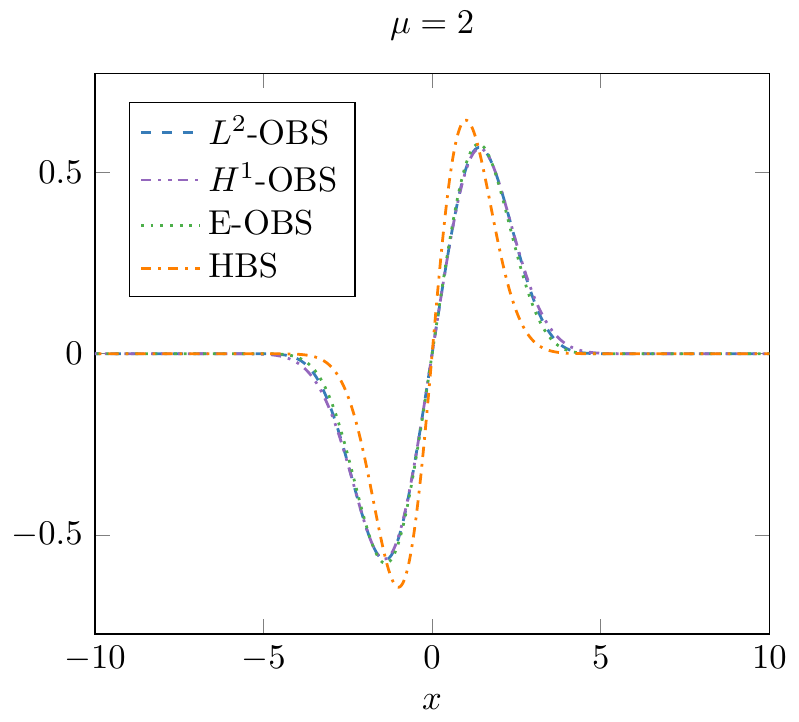}
    \\ \\
    $N_b = 3$ \\
    \includegraphics[width=0.25\linewidth]{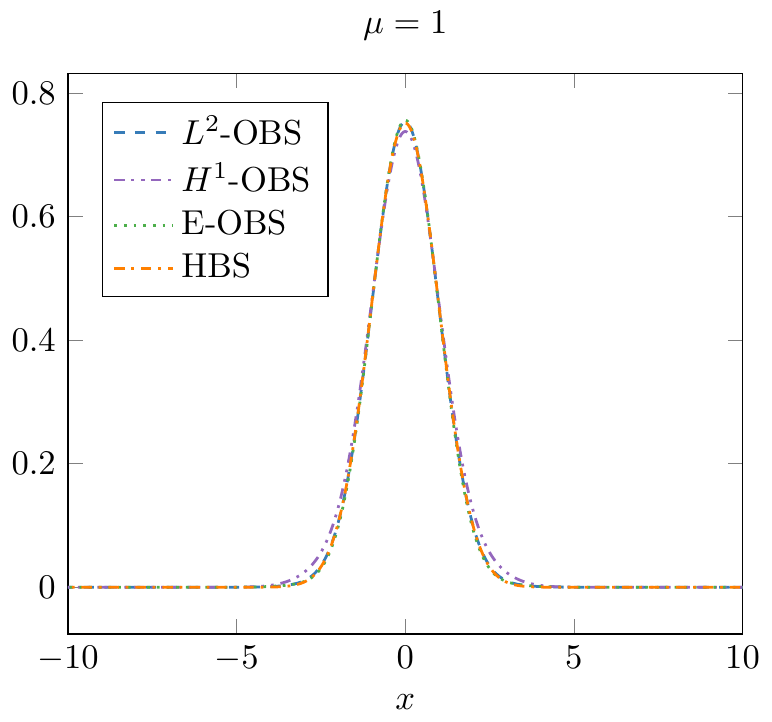}\includegraphics[width=0.25\linewidth]{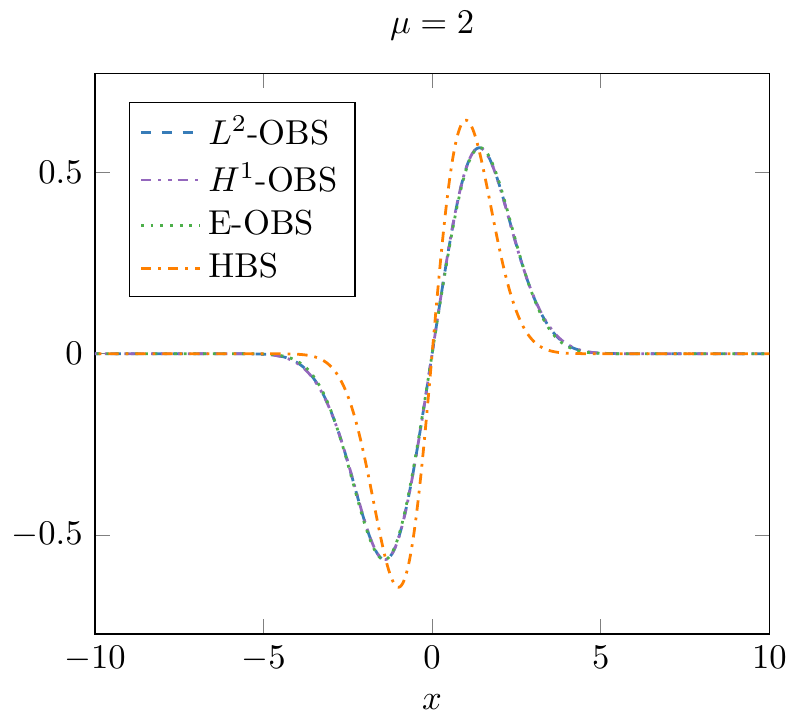}\includegraphics[width=0.25\linewidth]{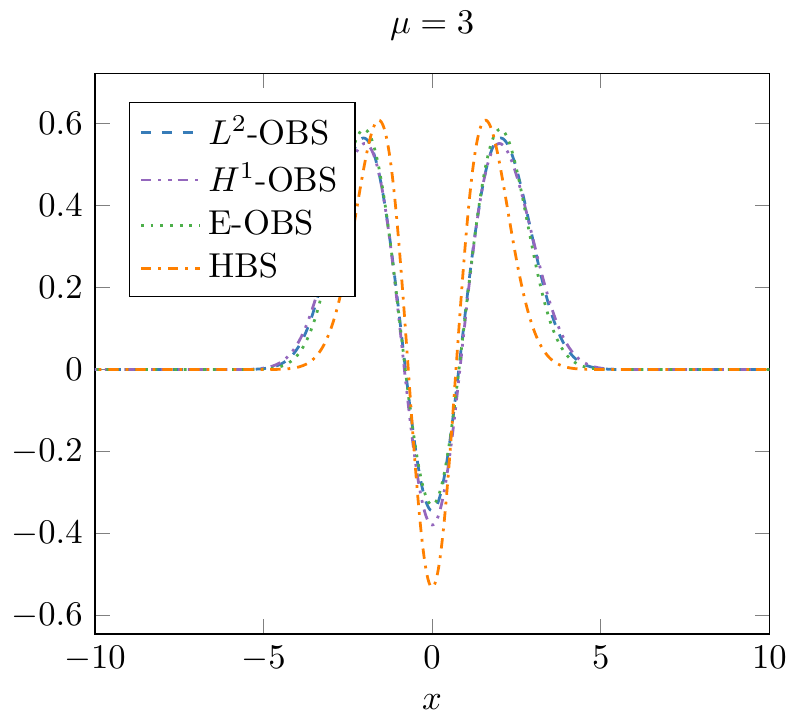}
    \\ \\
    $N_b = 4$ \\
    \includegraphics[width=0.25\linewidth]{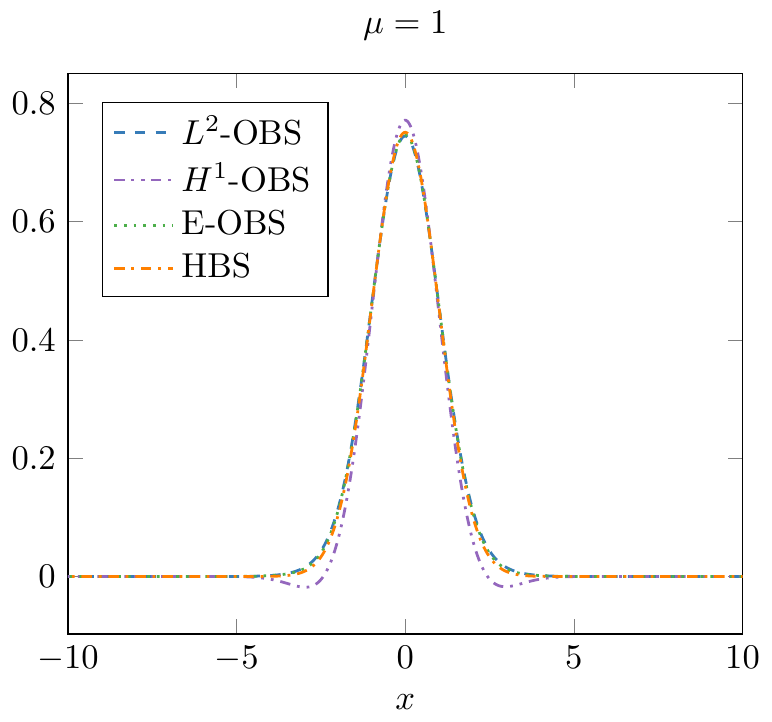}\includegraphics[width=0.25\linewidth]{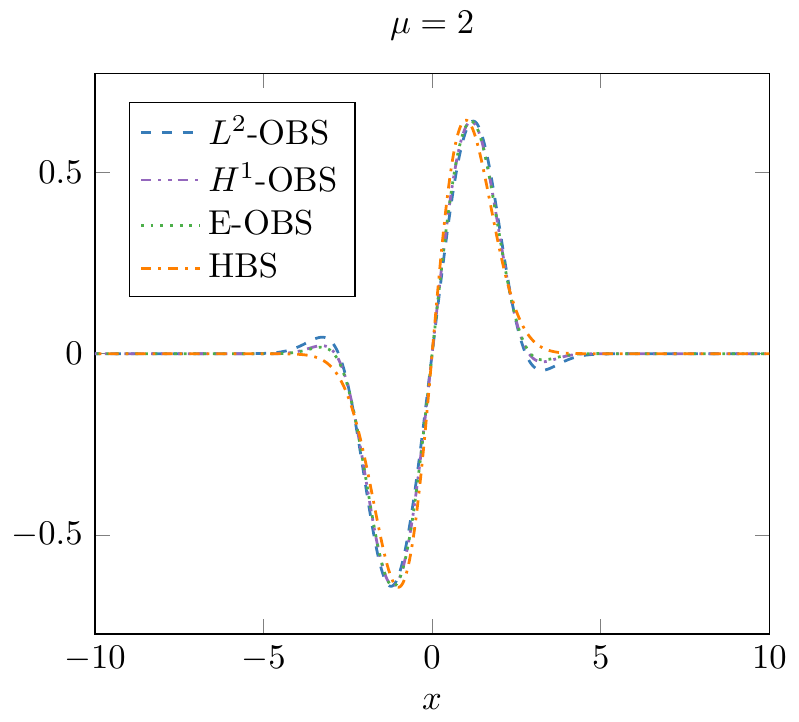}\includegraphics[width=0.25\linewidth]{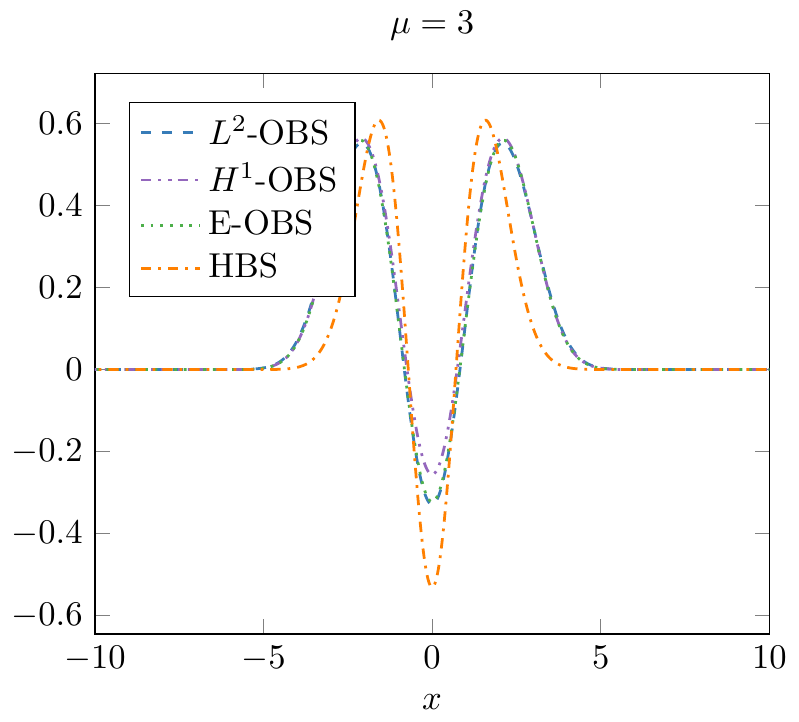}\includegraphics[width=0.25\linewidth]{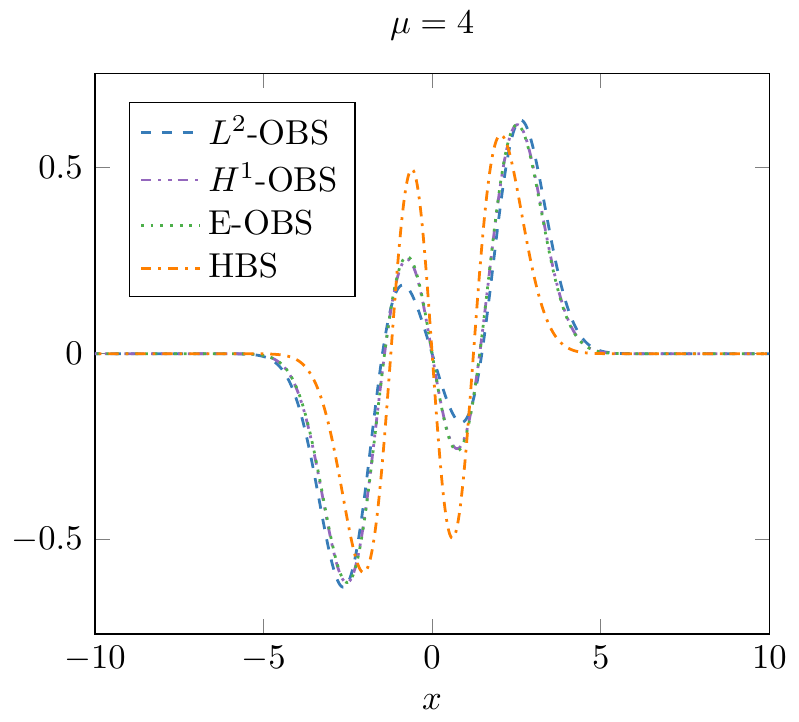}
    \\ \\
  \end{tabular}
  \caption{Optimal basis functions for different criteria, each of them being optimized for different values of $N_b$.}
  \label{fig:basis}
\end{figure}

Finally, for the sake of completeness, we plot in \autoref{fig:basis} the different basis functions built with each criterion for different values of $N_b$, confirming again the previous observations that the optimal basis functions are quite close to the standard Hermite basis functions.

\medskip

The main conclusion of these observations is that, for $N_b$ large enough, there is no real difference between the proposed criteria. Still, if the bases we built do not seem to be very different from the standard Hermite basis (\autoref{fig:basis}), building optimal bases allows to increase accuracy on the quantities of interest we focused on by one order of magnitude in average.

\vfill
\begin{table}[h!]
  \centering

  \begin{tabular}{@{}ccc@{}}
    Value of $J_{L^2}$ for the different basis sets & &
    Value of $J_{H^1}$ for the different basis sets \\ \\
    \begin{tabular}{@{}ccccc@{}}
      \toprule
      Basis     & $N_b=1$ & $N_b=2$ & $N_b=3$ & $N_b=4$ \\ \midrule
      HBS       & -7.40829 & -7.70051 & -7.74312 & -7.77138 \\
      $L^2$-OBS & -7.43954 & -7.76479 & -7.77725 & -7.77773 \\
      $H^1$-OBS & -7.43928 & -7.76466 & -7.77724 & -7.77772 \\
      E-OBS     & -7.39410 & -7.76425 & -7.77720 & -7.77772 \\ \bottomrule
    \end{tabular} & &
    \begin{tabular}{@{}ccccc@{}}
      \toprule
      Basis     & $N_b=1$ & $N_b=2$ & $N_b=3$ & $N_b=4$ \\ \midrule
      HBS       & -10.5613 & -11.0566 & -11.1451 & -11.2402 \\
      $L^2$-OBS & -10.6256 & -11.2338 & -11.2630 & -11.2650 \\
      $H^1$-OBS & -10.6265 & -11.2342 & -11.2630 & -11.2651 \\
      E-OBS     & -10.5334 & -11.2313 & -11.2626 & -11.2650 \\ \bottomrule
    \end{tabular}
  \end{tabular}

  \vspace{1cm}
  \begin{tabular}{@{}c@{}}
    Value of $J_{E}$ for the different basis sets\\ \\
    \begin{tabular}{@{}ccccc@{}}
      \toprule
      Basis     & $N_b=1$ & $N_b=2$ & $N_b=3$ & $N_b=4$ \\ \midrule
      HBS       & 3.77956$\times10^{-2}$ & 3.98301$\times10^{-3}$ & 1.86537$\times10^{-3}$ & 1.35309$\times10^{-4}$ \\
      $L^2$-OBS & 6.52016$\times10^{-2}$ & 2.18282$\times10^{-4}$ & 1.01365$\times10^{-6}$ & 3.22260$\times10^{-8}$ \\
      $H^1$-OBS & 6.83537$\times10^{-2}$ & 2.40548$\times10^{-4}$ & 1.27251$\times10^{-6}$ & 3.91885$\times10^{-8}$ \\
      E-OBS     & 3.69610$\times10^{-2}$ & 1.92087$\times10^{-4}$ & 6.93394$\times10^{-7}$ & 2.54014$\times10^{-8}$ \\ \bottomrule
    \end{tabular}
  \end{tabular}

  \vspace{1cm}
  \begin{tabular}{@{}c@{}}
    L--BFGS iterations\\ \\
    \begin{tabular}{@{}ccccc@{}}
      \toprule
      Basis       & $N_b=1$ & $N_b=2$ & $N_b=3$ & $N_b=4$ \\ \midrule
      $L^2$-OBS & 4       & 13      & 48      & 219     \\
      $H^1$-OBS & 7       & 17      & 235     & not converged after 500 it   \\
      E-OBS     & 6       & 19      & 52      & 134     \\ \bottomrule
    \end{tabular}
  \end{tabular}
  \vspace{1cm}
  \caption{(Top \& Middle) Values of the different criteria for the HBS and optimal bases, for increasing values of $N_b$. (Bottom) Number of iterations of L--BFGS required for each criterion to achieve convergence up to requested tolerance ($10^{-7}$ on the $\ell^2$-norm of the gradient).}
  \label{tab:criteria}
\end{table}
\vfill

\clearpage

\subsection{Influence of numerical parameters}

\subsubsection{Random starting points}

In \autoref{subsec:num_res}, we used the first $N_b$ Hermite functions as a starting point for the optimization procedures. We obtain the same solutions if we start from a random matrix $R$ on the Stiefel manifold, in the sense that the optimal values reached for each criterion are the same, as well as the error plots. However, the L--BFGS algorithm requires more iterations to converge. The basis functions obtained from the optimization algorithms are different from those observed in \autoref{fig:basis}, but still span the same space as the variational solutions are equal.

\subsubsection{Extrapolating the parameter space $\Ic$}
\label{subsubsec:extrapolation}

In \autoref{subsec:num_res}, we chose a probability measure $\PP$ supported in the interval $[1.5,5]$ in order to avoid conditioning issues.
Indeed, taking smaller values of $a$ results in the L--BFGS algorithm having convergence problems when $N_b$ increases. This phenomenon was observed already for $N_b=3$ or $N_b=4$ when including $a=1$ in the support of $\PP$. In practice, this problem can be solved by using preconditioning or getting rid of overcompleteness by pre-processing the basis $\chi_a$ (e.g. selecting a smaller basis by filtering out the very small singular values of the original overlap matrix), but for brevity we will not elaborate further in this direction.

\medskip

\begin{figure}[h!]
  \centering
  \includegraphics[width=0.9\linewidth]{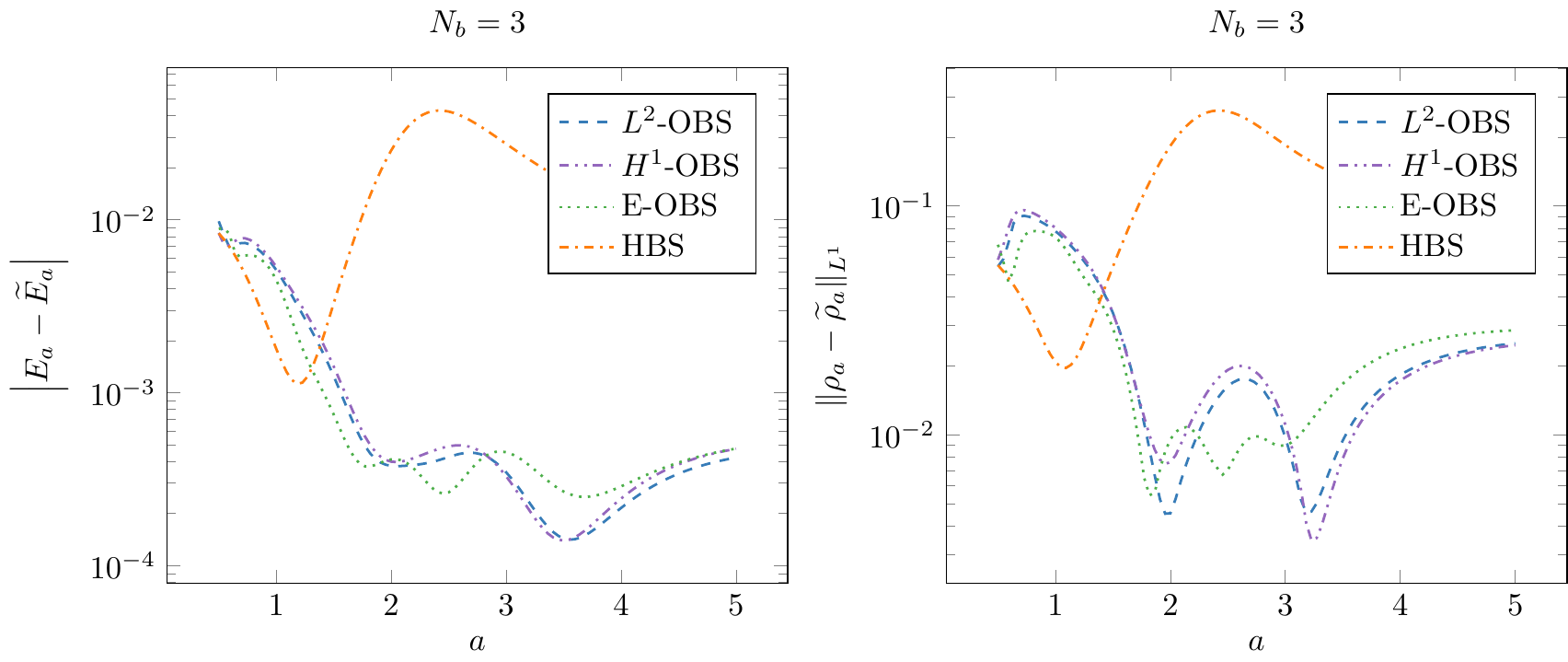}
  \caption{Energy and densities error with extrapolation up to $a=0.5$, with basis functions optimized on $\Ic=[1.5,5]$}.
  \label{fig:extrapolation}
\end{figure}

However, once we have computed optimal bases for a reasonable interval $\Ic$, it is possible to use these bases to solve the variational problem \eqref{eq:pb_init_galerkin} and extrapolate the energy and the density to smaller values of $a$ that are not in the set $\Ic$. The results are plotted in \autoref{fig:extrapolation}. We notice that the quantities of interest are better approximated on $\Ic = [1.5,5]$, but for smaller $a$'s, there is no more gain in accuracy with respect to the standard HBS.

\subsubsection{Choice of the probability $\PP$}
\label{subsubsec:choiceP}

The major drawback of our AO basis optimization lies in the necessity to compute very accurate reference solutions for all configurations in the support of $\PP$. This is not an issue for our 1D toy model but it can be very time consuming for real systems if the support of $\PP$ is too large. It is therefore crucial to reduce as much as possible the support of $\PP$.

\medskip

\begin{figure}[h!]
  \centering
  \includegraphics[width=\linewidth]{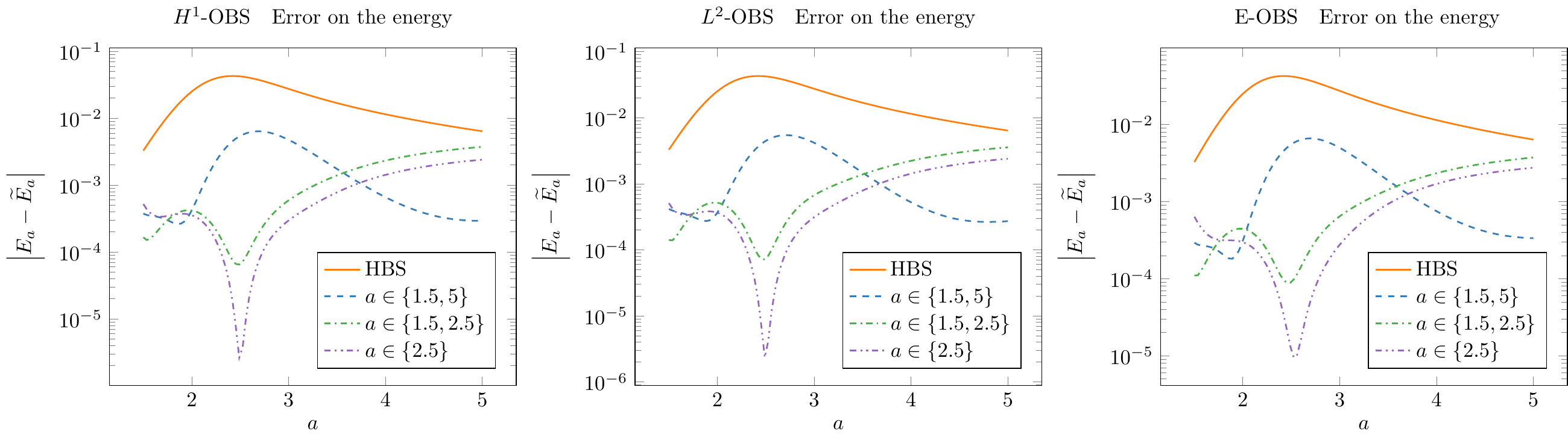}
  \includegraphics[width=\linewidth]{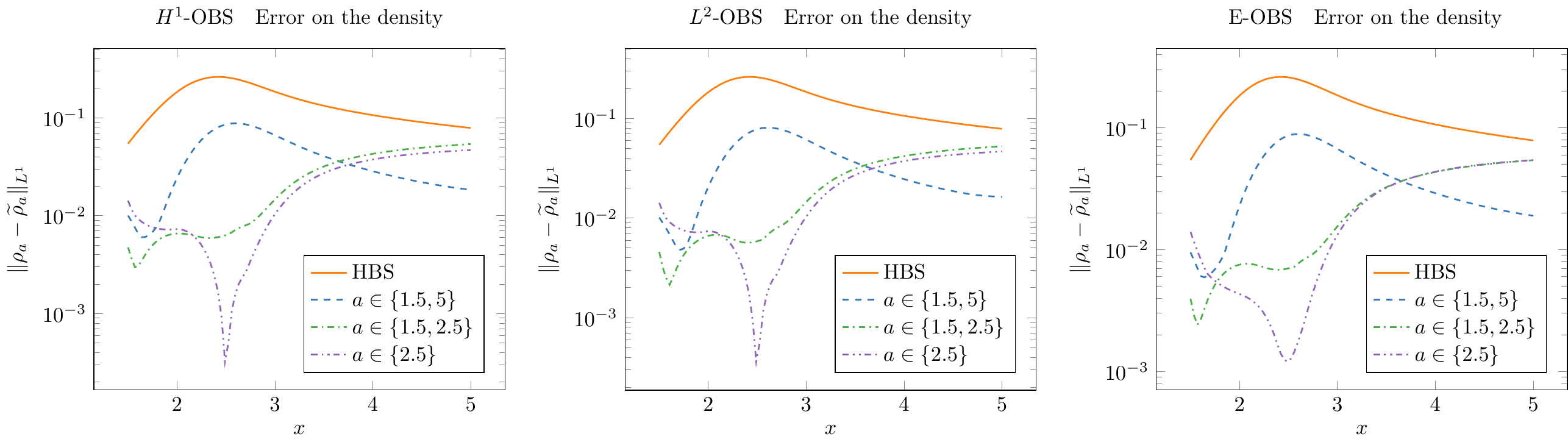}
  \caption{Error plots for probability measures $\PP$ corresponding to very sparse samplings of the interval $\Ic=[1.5,5]$: i) the two endpoints of $\Ic$ ii) two points near the equilibrium distance and iii) one point near the equilibrium distance.
    (Top line) Error on energy. (Bottom line) Error on density in $L^1$ norm.
    (Left) OBS for $J_{H^1}$. (Middle) OBS for $J_{L^2}$. (Right) OBS for $J_E$.
    The \enquote{$a$} in legends are the sampled configurations $a$.}
  \label{fig:compare_sampling}
\end{figure}

In this section, we study the influence of the probability measure $\PP$ on the quality of the optimized bases. For simplicity, we restrict ourselves to uniform samplings of the interval $\Ic=[1.5,5]$. Numerical tests show that increasing the sample size above the reference sampling with $N_{\rm c}=10$ points used in \autoref{subsec:num_res} (see Eq.~\eqref{eq:probaP}) brings no significant accuracy improvement. Therefore we chose to investigate in the following the performance of the optimal AO basis sets obtained with very sparse sampling. \autoref{fig:compare_sampling} pictures the error of approximation of the dissociation curve and densities for three samplings: first, the two extreme points of the interval $\Ic=[1.5,5]$; second, two points around the equilibrium distance $a_0\simeq 1.925$ ; third, a single point near the equilibrium distance. All curves are plotted for a fixed number of basis functions $N_b=3$.

\medskip

It appears that the latter sampling already provides satisfactory accuracy. The criteria $J_{L^2}$ and $J_{H^1}$ are equal to $-5\times 10^{-6}$ for optimized basis to be compared with $-1.8\times 10^{-3}$ for standard HBS. Hence they provide a gain of accuracy in energy of three orders of magnitude over the whole dissociation curve.

\subsubsection{Number of Hilbert basis functions}
We now take the same setting as in \autoref{subsec:num_res}, except that we
set $\Nc = 5$ instead of $\Nc = 10$. This provides similar results as those
collected in \autoref{tab:criteria}, see \autoref{tab:criteria_N5}. However,
the values of the criteria $J_A$ and $J_E$ are higher than for $\Nc = 10$, in
particular for $N_b=4$, where criterion $J_A$ cannot be optimized further than
$-10^{-5}$, which makes sense as the space over which the optimization
algorithms are performed is smaller. Calculations with $\Nc = 15$ were also
performed: for $N_b=1,2,3$, the criteria are slightly
improved but for $N_b=4$, convergence issues were noticed, due to ill
conditioning of the overlap matrices for $a=1.5$ as the number $\Nc$ of
functions used to describe the optimal bases is larger.

\begin{table}[h!]

  \begin{tabular}{@{}ccc@{}}
    Value of $J_{L^2}$ for the different basis sets & &
    Value of $J_{H^1}$ for the different basis sets \\ \\
    \begin{tabular}{@{}ccccc@{}}
      \toprule
      Basis     & $N_b=1$ & $N_b=2$ & $N_b=3$ & $N_b=4$ \\ \midrule
      HBS       & -7.40829 & -7.70051 & -7.74312 & -7.77138 \\
      $L^2$-OBS & -7.43933 & -7.76304 & -7.77554 & -7.77618 \\
      $H^1$-OBS & -7.43923 & -7.76258 & -7.77525 & -7.77612 \\
      E-OBS     & -7.39401 & -7.76259 & -7.77545 & -7.77615 \\ \bottomrule
    \end{tabular} & &
    \begin{tabular}{@{}ccccc@{}}
      \toprule
      Basis     & $N_b=1$ & $N_b=2$ & $N_b=3$ & $N_b=4$ \\ \midrule
      HBS       & -10.5613 & -11.0566 & -11.1451 & -11.2402  \\
      $L^2$-OBS & -10.6237 & -11.2225 & -11.2541 & -11.2577 \\
      $H^1$-OBS & -10.6240 & -11.2244 & -11.2555 & -11.2581 \\
      E-OBS     & -10.5328 & -11.2234 & -11.2547 & -11.2580 \\ \bottomrule
    \end{tabular}
  \end{tabular}

  \vspace{1cm}
  \begin{tabular}{@{}c@{}}
    Value of $J_{E}$ for the different basis sets\\ \\
    \begin{tabular}{@{}ccccc@{}}
      \toprule
      Basis     & $N_b=1$ & $N_b=2$ & $N_b=3$ & $N_b=4$ \\ \midrule
      HBS       & 3.77956$\times10^{-2}$ & 3.98301$\times10^{-3}$ & 1.86537$\times10^{-3}$ & 1.35309$\times10^{-4}$ \\
      $L^2$-OBS & 6.43832$\times10^{-2}$ & 2.46466$\times10^{-4}$ & 1.58667$\times10^{-5}$ & 1.01128$\times10^{-5}$ \\
      $H^1$-OBS & 6.13025$\times10^{-2}$ & 2.45930$\times10^{-4}$ & 1.62235$\times10^{-5}$ & 1.00611$\times10^{-5}$ \\
      E-OBS     & 3.69681$\times10^{-2}$ & 1.30365$\times10^{-4}$ & 1.41935$\times10^{-5}$ & 9.74560$\times10^{-6}$ \\ \bottomrule
    \end{tabular}
  \end{tabular}
  \vspace{1cm}
  \caption{Value of the different criteria for the different local (optimized and Hermite) bases, with $\Nc=5$ and increasing values of $N_b$.}
  \label{tab:criteria_N5}
\end{table}

\section*{Acknowledgements}
The authors thank Etienne Polack for fruitful discussions.
This project has received funding from the
European Research Council (ERC) under the European Union's Horizon 2020
research and innovation programme (grant agreement EMC2 No 810367). This work
was supported by the French ‘Investissements d’Avenir’ program, project Agence
Nationale de la Recherche (ISITE-BFC) (contract ANR-15-IDEX-0003).

\clearpage
\appendix
\section*{Appendix}

In this appendix, we will use extensively the two symmetries of the trace: for any matrices $M$ and $N$ such that $MN$ and $NM$ are defined,
\[
  {\rm Tr}(MN) = {\rm Tr}(NM)\quad\mbox{and}\quad {\rm Tr}(M^T) = {\rm Tr}(M).
\]

\subsubsection*{Computation of the gradient of \texorpdfstring{$J_A$}{}}
Let $R,H \in \RR^{\Nc\times N_{b}}$ and define $I_H=\begin{pmatrix}
  H & 0 \\ 0 & H
\end{pmatrix}$.
One has
\begin{equation}
  \label{eq:grad_K_1}
  \begin{split}
    J_A(R+H) - J_A(R) & = - \int_\Omega {\rm Tr}\left(M_A^{\rm offline}(a)\left(2I_R [S^A (B_a I_R^{\,})]^{-1} I_H^T + I_R^{\,}\,\left[{\rm d}[S^A]^{-1} (B_a I_R) \cdot (B_a
          I_H^{\,})\right] I_R^T\right)\right) \, \d\PP(a)\\
    & + O(\Vert H \Vert^2)
  \end{split}
\end{equation}
Considering that
\[
  (M+H)^{-1} - M^{-1} = -M^{-1}HM^{-1} + O(\Vert H \Vert^2)\mbox{ and } S^A(BI_{R+H}) - S^A(B_aI_R) = I_H^TS^A(B)I_R^{\,} + I_R^T S^A(B) I_H^{\,} + O(\Vert H \Vert^2),
\]
it follows from the chain rule that
\begin{align*}
  {\rm d}[S^A]^{-1}(B_aI_R)\cdot(B_aI_H) &= - [S^A(B_aI_R)]^{-1}\left(I_H^T S^A(B_a)I_R^{\,} + I_R^TS^A(B_a) I_H^{\,} \right)[S^A(B_aI_R)]^{-1}.
\end{align*}
From this computation, we obtain that the integrand in expression \eqref{eq:grad_K_1} writes for all $a$
\begin{align}\label{eq:grad_J_A_intermediate}
  &2{\rm Tr}\left(M_A^{\rm offline}(a)\left[I_R [S^A(B_aI_R^{\,})]^{-1}I_H^T - I_R [S^A(B_aI_R)]^{-1}I_H^T S^A(B_a)I_R^{\,}[S^A(B_aI_R)]^{-1}I_R^T \right] \right)\nonumber\\[0.15cm]
  &= 2{\rm Tr}\left(M_A^{\rm offline}(a)I_R [S^A(B_aI_R^{\,})]^{-1}I_H^T - I_H^T S^A(B_a)I_R^{\,}[S^A(B_aI_R)]^{-1}I_R^TM_A^{\rm offline}(a)I_R [S^A(B_aI_R)]^{-1} \right).
\end{align}

The idea is now to write the expression \eqref{eq:grad_J_A_intermediate} as the inner product of $H$ with a given matrix of $\mathbb{R}^{\Nc\times N_b}$, which we will identify as the integrand of the gradient of $K_I$. Changing from $I_H$ to $H$ imposes to decompose each matrix by block and to write the trace in \eqref{eq:grad_J_A_intermediate} as the sum of traces over the diagonal blocks. To this end we introduce the superscripts "$++$", "$+-$", "$-+$" and "$--$" associated with one of the four identically shaped blocks of a generic matrix
\begin{equation}
  \label{eq:pp_mm}
  M=\begin{pmatrix}M^{++}&M^{+-}\\ M^{-+}&M^{--}\end{pmatrix}.
\end{equation}
Expression \eqref{eq:grad_J_A_intermediate} therefore immediately reads
\begin{align}
  & 2{\rm Tr}\biggl(I_H^T\underbrace{\left[M_A^{\rm offline}(a)I_R^{\,}[S^A(B_aI_R^{\,})]^{-1}- S^A(B_a)I_R^{\,}[S^A(B_aI_R)]^{-1}I_R^T M_A^{\rm offline}(a)I_R[S^A(B_aI_R)]^{-1} \right]}_{\displaystyle M_A(a,R)}\biggr)\nonumber\\[0.15cm]
  & \hspace{2cm} = 2{\rm Tr}\left(H^T \left(M_A(a,R)^{++} + M_A(a,R)^{--} \right)\right).
\end{align}
One can verify that $M_A(a,R)^{++} +  M_A(a,R)^{--}$ is in $\mathbb{R}^{\Nc \times N_b}$ and we conclude by identification that
\begin{equation}
  \boxed{\nabla J_A(R) = -2\int_\Omega \left(M_A(a,R)^{++} + M_A(a,R)^{--}\right){\rm d}\PP(a).}
\end{equation}

\subsubsection*{Computation of the gradient of \texorpdfstring{$J_E$}{}}

Let $R,H \in \RR^{\Nc\times N_b}$ and define $I_H=\begin{pmatrix}
  H & 0 \\ 0 & H
\end{pmatrix}$.
We immediately have that
\begin{equation}
  \label{eq:gradient_J_E}
  \nabla J_E(R) = -2\int_\Omega \nabla E_a(R)\left(E_a^{\,} - E_a(R)\right){\rm d} \PP(a),
\end{equation}
where
\begin{equation}
  E_a(R) = {\rm Tr}\left(C_a(R) (C_a(R))^T\times \Hc_a(R) \right),
\end{equation}
with $C_a(R)$ defined in \autoref{subsec:discrete_setting} and
$\Hc_a(R) \coloneqq I_R^T M_E^{\rm offline}(a)I_R^{\,}$. Therefore, if we define $\Ec_a(R,C) = {\rm Tr}(CC^T\Hc_a(R))$, then
$E_a(R) = \Ec_a(R,C_a(R))$ and we have, by the chain rule,
\begin{equation*}
  \nabla E_a(R)\cdot H = \nabla_R \Ec_a(R,C_a(R))\cdot H +
  \nabla_C \Ec_a(R,C_a(R))\cdot\left( {\rm d}  C_a(R)\cdot H\right).
\end{equation*}
We now detail the computations of the two gradients of $\Ec_a$, namely $\nabla_R \Ec_a$ and $\nabla_C \Ec_a$.

\paragraph{\underline{Computation of the first gradient $\nabla_R \Ec_a$}}
Using notation \eqref{eq:pp_mm}, we introduce
\begin{equation*}
  M_a \coloneqq M_E^{\rm offline}(a) \text{ and }
  \Sigma(H) \coloneqq I_H^T M_a I_R =
  \begin{pmatrix}
    H^TM_a^{++}R & H^TM_a^{+-}R\\
    H^{T^{\phantom{T}}}M_a^{-+}R & H^TM_a^{--}R\\
  \end{pmatrix} \in \RR^{(2N_b)\times (2N_b)},
\end{equation*}
so that, with $P = CC^T$,
\begin{equation*}
  \begin{split}
    {\rm Tr}\left( P[{\rm d}\Hc_a(R)\cdot H] \right) &=
    {\rm Tr}\left(P[\Sigma(H) + \Sigma(H)^T]\right) =
    2{\rm Tr}(P\Sigma(H)) \\ &=
    2{\rm Tr}\left(H^T\left(
        M_a^{++}RP^{++} + M_a^{-+}RP^{+-}
        + M_a^{+-}RP^{-+} + M_a^{--}RP^{--}
      \right)\right).
  \end{split}
\end{equation*}
In the end,
\begin{equation*}
  {
    \nabla_R\Ec_a(R,C) =
    2\left(M_a^{++}R(CC^T)^{++} + M_a^{+-}R(CC^T)^{-+}
      + M_a^{-+}R(CC^T)^{+-} + M_a^{--}R(CC^T)^{--}\right) \in \RR^{\Nc \times N_b}.
  }
\end{equation*}

\paragraph{\underline{Computation of the second gradient $\nabla_C \Ec_a$}}
The Euler--Lagrange equation of the minimization problem \eqref{eq:E_crit_JE} yields that there
exist a symmetric matrix $\Lambda_a(R)\in\RR^{2\times2}$ such that
\begin{equation*}
  \nabla_C \Ec_a(R,C_a(R)) = 2 \Hc_a(R) = 2S(B_aI_R)C_a(R)\Lambda_a(R),
\end{equation*}
where $\Lambda_a(R)$ is actually a diagonal matrix whose diagonal is composed
of the two lowest eigenvalues of $\Hc_a(R)$.  Moreover, if we differentiate the
constraint $C_a(R)^TS(B_a I_R)C_a(R) = \Id_2$, we get
\begin{equation*}
  C_a(R)^TS(B_aI_R)({\rm d} C_a(R)\cdot H) +
  \left({\rm d} C_a(R)\cdot H\right)^TS(B_aI_R)C_a(R) = -
  C_a(R)^T(\d S(B_aI_R)\cdot H)C_a(R),
\end{equation*}
so that
\begin{equation*}
  \begin{split}
    \nabla_C \Ec_a(R,C_a(R))\cdot\left( {\rm d} C_a(R)\cdot H\right) &=
    2{\rm Tr}\left( (S(B_aI_R)C_a(R)\Lambda_a(R))^T ({\rm d} C_a(R)\cdot H)
    \right) \\ &= -{\rm Tr}\left(
      \left({\rm d} S(B_aI_R)\cdot H\right)C_a(R)\Lambda_a(R)C_a(R)^T
    \right).
  \end{split}
\end{equation*}
Now, let us recall that
\begin{equation*}
  {\rm d} S(B_aI_R)\cdot H = I_H^TS(B_a)I_R + I_R^TS(B_a)I_H.
\end{equation*}
Thus, by denoting $Q_a(R) = C_a(R)\Lambda_a(R)C_a(R)^T$, we get that
\begin{equation*}
  \begin{split}
    \nabla_C \Ec_a(R,C_a(R))\cdot\left( {\rm d} C_a(R)\cdot H\right)
    =&- 2
    {\rm Tr}\left(
      H^T\left( S(B_a)^{++}RQ_a(R)^{++} + S(B_a)^{+-}RQ_a(R)^{-+}\right.\right.\\
    &\phantom{\rm Tr}+\left.\left.S(B_a)^{-+}RQ_a(R)^{+-} + S(B_a)^{--}RQ_a(R)^{--} \right)
    \right)
  \end{split}
\end{equation*}
which ends the computations of the second gradient.

\paragraph{\underline{Final gradient}}

Compiling the computations of the two previous paragraphs, we obtain
\begin{equation}
  \boxed{
    \begin{split}
      \nabla_R E_a(R) &=
      2\left(M_a^{++}RP_a(R)^{++} + M_a^{+-}RP_a(R)^{-+}
        + M_a^{-+}RP_a(R)^{+-} + M_a^{--}RP_a(R)^{--}\right) \\
      &\phantom{=}-2 \left(S_a^{++}RQ_a(R)^{++} + S_a^{+-}RQ_a(R)^{-+}
        + S_a ^{-+}RQ_a(R)^{+-} + S_a^{--}RQ_a(R)^{--}\right)
    \end{split}
  }
\end{equation}
where $P_a(R)=C_a(R)C_a(R)^T$, $M_a = M_A^{\rm offline}(a)$, $S_a = S(B_a)$ and $Q_a(R) = C_a(R) \Lambda_a(R) C_a(R)^T$, and the gradient of $J_E$ is computed with $\eqref{eq:gradient_J_E}$.

\bibliographystyle{siam}
\bibliography{biblio.bib}

\end{document}